\documentclass[12pt]{amsart}
\usepackage{amssymb,amsmath,amsthm}
\numberwithin{equation}{section}
\begin{document}

\theoremstyle{plain}
\newtheorem{theorem}{Theorem}[section]
\newtheorem{lemma}[theorem]{Lemma}
\newtheorem{proposition}[theorem]{Proposition}
\newtheorem{corollary}[theorem]{Corollary}
\newtheorem{conjecture}[theorem]{Conjecture}

\theoremstyle{definition}
\newtheorem*{definition}{Definition}

\theoremstyle{remark}
\newtheorem*{remark}{Remark}
\newtheorem{example}{Example}[section]
\newtheorem*{remarks}{Remarks}

\newcommand{\cc}{{\mathbb C}}
\newcommand{\qq}{{\mathbb Q}}
\newcommand{\rr}{{\mathbb R}}
\newcommand{\nn}{{\mathbb N}}
\newcommand{\zz}{{\mathbb Z}}
\newcommand{\pp}{{\mathbb P}}
\newcommand{\al}{\alpha}
\newcommand{\be}{\beta}
\newcommand{\ga}{\gamma}
\newcommand{\ze}{\zeta}
\newcommand{\om}{\omega}
\newcommand{\ep}{\epsilon}
\newcommand{\la}{\lambda}
\newcommand{\de}{\delta}
\newcommand{\De}{\Delta}
\newcommand{\si}{\sigma}
\newcommand{\ta}{\theta}

\title{Some new canonical forms for polynomials}

\author{Bruce Reznick}
\address{Department of Mathematics and Center for Advanced Study, University of 
Illinois at Urbana-Champaign, Urbana, IL 61801} 
\email{reznick@math.uiuc.edu}
 \date{\today}
\begin{abstract}
We give some new canonical representations for forms over $\cc$. 
For example, a general binary quartic form can be written as the  square of
a quadratic form plus the fourth power of a linear form. A general cubic form in
$(x_1,\dots,x_n)$ can be written uniquely as a sum of the cubes of
linear forms $\ell_{ij}(x_i,\dots,x_j)$, $1 \le i \le j \le n$. A
general ternary quartic form is the sum of the square of a quadratic form and
three fourth powers of linear forms. The methods are classical and elementary.
\end{abstract}

\maketitle
\section{Introduction and Overview}
\subsection{Introduction}
Let $H_d(\cc^n)$ denote the $N(n,d) = \binom{n+d-1}d$-dimensional
vector space of complex forms of degree $d$ in $n$ variables, or
{\it $n$-ary $d$-ic forms}. One of
the   major accomplishments of 19th century algebra was the
discovery of canonical forms for certain classes of $n$-ary $d$-ics,
especially as the sum of $d$-th power of linear forms. By a {\it
  canonical form} we mean
a polynomial $F(t;x)$ in two sets of variables, $t \in \cc^{N(n,d)}$ and $x \in
\cc^n$, with the property that for general $p \in H_d(\cc^n)$,
there exists $t$ so that $p(x) = F(t;x)$. Put another way, the set
$\{F(t;x): t \in \cc^{N(n,d)}\}$ is a Zariski open set in
$H_d(\cc^n)$.

In this paper, we present some new canonical
forms, whose main novelty  is that they involve intermediate powers 
of forms of  higher degree, or forms with a restricted set of monomials.
(These variations have been suggested by Hilbert's study of ternary
quartics \cite{H}, which led to his 17th problem,  as well as by a
remarkable  theorem of B. Reichstein \cite{Re1}  on cubic forms.) These
expressions, are less
susceptible to apolarity arguments than  
the traditional canonical forms, and lead naturally to (mostly open) enumeration
questions. 

To take a simple, yet familiar example, 
\begin{equation}\label{E:toy}
F(t_1,t_2,t_3;x,y) = (t_1x + t_2 y)^2 + (t_3 y)^2
\end{equation}
is a canonical form for binary quadratic forms. By the usual completion of
squares, $p(x,y) = ax^2 + 2b xy + cy^2$ can be put into \eqref{E:toy}  
for $t_1 =
\sqrt a$, $t_2 =  b/t_1$ and $t_3^2 = c - t_2^2$. Many of the 
examples in this paper can be viewed as imperfect 
attempts to generalize \eqref{E:toy}.

In 1851, Sylvester \cite{S1,S2} presented a family of canonical forms for
binary forms in all degrees.

\begin{theorem}[Sylvester's Theorem]\label{Sylvcan} 
\ 

\noindent (i) A general binary form $p$ of odd degree $2s-1$ can be
written as
\begin{equation}\label{E:Sylvcanodd}
p(x,y) = \sum_{j=1}^s (\al_j x + \be_j y)^{2s-1}.
\end{equation}
(ii) A general binary form $p$ of even degree $2s$ can be written as
\begin{equation}\label{E:Sylvcaneven}
p(x,y) = \la x^{2s} + \sum_{j=1}^s (\al_j x + \be_j y)^{2s}.
\end{equation}
for some $\la \in \cc$.
\end{theorem}
The somewhat unsatisfactory nature of the asymmetric summand in
\eqref{E:Sylvcaneven} has been the inspiration for other canonical
forms for binary forms of even degree.

Another familiar canonical form  is the generalization of \eqref{E:toy} into the
upper-triangular expression for quadratic forms, found by
repeated completion of the square:
\begin{theorem}\label{uppertri}
A general quadratic form $p \in H_2(\cc^n)$ can be written as:
\begin{equation}\label{E:upper}
  p(x_1,\dots,x_n) =\sum_{k=1}^n (t_{k,k}x_k + t_{k,k+1}x_{k+1}
  + \dots + t_{k,n}x_n)^2, \quad t_{k,\ell} \in \cc. 
\end{equation}
\end{theorem}
The expression in \eqref{E:upper} is unique, up to the signs of
the linear forms.

There are two ways to verify that a candidate expression $F(t;x)$ is,
in fact, a canonical
form. One is the classical non-constructive method based on the
existence of a point at which the Jacobian matrix has full rank. (See Corollary 
\ref{jake}, and see Theorem \ref{LW} for the apolar version.)
 Lasker \cite{Lasker} attributes the underlying idea 
 to Kronecker and L\"uroth -- see \cite[p.208]{W}.   

Ideally, however, a canonical form can be derived constructively, and the
number of different representations can thereby be determined.  The convention
in this paper will be that two representations  are the same if they
are equal, up 
to a permutation of like summands and with the identification of $f^k$
and $(\zeta f)^k$ when $\zeta^k = 1$. The representation  in
\eqref{E:Sylvcanodd} 
is unique in this sense, even though there are $s!\cdot(2s-1)^s$
different $2s$-tuples $(\al_1,\be_1,\dots,\al_s,\be_s)$ for which
\eqref{E:Sylvcanodd} is valid. 

In addition to Theorem \ref{Sylvcan}, another motivational example for
this paper is a
remarkable canonical form for cubic forms found by Reichstein \cite{Re1} 
in 1987, which can be thought of as a ``completion of the cube''. 
\begin{theorem}[Reichstein]\label{reichcan}
A general cubic $p\in H_3(\cc^n)$ can be written uniquely as 
\begin{equation}\label{E:reich}
 p(x_1,\dots,x_n) = \sum_{k=1}^n \ell_k^3(x_1,\dots,x_n) + q(x_3,\dots,x_n), 
\end{equation}
where $\ell_k \in H_1(\cc^n)$ and  $q \in H_3(\cc^{n-2})$.
\end{theorem}
This is a canonical form, provided $q$ is viewed as a $t$-linear combination
of the monomials in $(x_3,\dots,x_n)$; since $N(n,3) = n^2 +
N(n-2,3)$, the constant 
count is right. Iteration (see \eqref{E:fullreich}) gives $p$ as a sum of
roughly $n^2/4$ cubes. The minimum from constant-counting, which is
justified by the Alexander-Hirschowitz Theorem \cite{AH}, is roughly $n^2/6$.
We give Reichstein's constructive proof of Theorem
\ref{reichcan} in section six. 

Here are some representative examples of the new canonical forms in
this paper. 

\begin{theorem}\label{slinkycan}
A general cubic form $p \in H_3(\cc^n)$ has a unique representation 
\begin{equation}\label{E:slinky}
  p(x_1,\dots,x_n) = \sum_{1 \le i \le j \le n} (t_{\{i,j\},i}x_i + \cdots
+ t_{\{i,j\},j}x_j)^3,
\end{equation}
where $t_{\{i,j\},k} \in \cc$.
\end{theorem}

\begin{theorem}\label{sextican}
A general binary sextic $p\in H_6(\cc^2)$ can be written as $p(x,y) = f^2(x,y) +
g^3(x,y)$, where $f\in H_3(\cc^2)$ is a  cubic form and $g \in
H_2(\cc^2)$ is a quadratic form. 
\end{theorem}

Theorem \ref{slinkycan}  has a constructive proof. Theorem
\ref{sextican} is in fact, a very special case of much deeper recent
results of V\'arilly-Alvarado. (See \cite{V1}, especially Theorem 1.2
and Remark 4.5, and Section 1.2 of \cite{V2}.) We  include it
because our proof, in the next section, is very short. 

Theorems \ref{Sylvcan} and \ref{sextican} are both special cases of a more
general class of canonical forms for $H_d(\cc^2)$, which is a
corollary of \cite[Theorem\ 4.4]{ER} (see Theorem \ref{ER}), but not
worked out explicitly there. 
\begin{theorem}\label{omnibus}
Suppose $d \ge 1$, $\{\ell_j: 1 \le j \le m\}$ is a fixed 
set of pairwise non-proportional linear forms, and
suppose $e_k \ | \ d$, $d > e_1 \ge \cdots \ge e_r$, $1 \le k \le r$, and 
\begin{equation}\label{E:book}
m + \sum_{k=1}^r(e_k+1) = d+1.
\end{equation}
 Then a general binary $d$-ic form $p \in H_d(\cc^2)$ can be written as
\begin{equation}\label{E:reveal}
p(x,y) = \sum_{j=1}^m t_j \ell_j^d(x,y) + \sum_{k=1}^r
f_k^{d/e_k}(x,y),\quad 
\end{equation}
where $t_j \in \cc$ and $\deg f_k = e_k$.
\end{theorem}
The condition $e_k < d$ excludes the vacuous case $m=0, r=1, e_1 = d$.
If each $e_k = 1$ and $r = \lfloor \frac {d+1}2
\rfloor$, then $m = d+1 - 2 \lfloor \frac {d+1}2 \rfloor \in \{0,1\}$ and
Theorem \ref{omnibus} becomes Theorem \ref{Sylvcan}; Theorem 
\ref{sextican} is Theorem  \ref{omnibus} in the special case 
$d=6, m=0, r=2, e_1=3, e_2=2$. As an example of a canonical form which
is unlikely to find a constructive proof: for a general   
$p \in H_{84}(\cc^2)$, there exist $f \in H_{42}(\cc^2),g \in
H_{28}(\cc^2)$ and  $h \in  H_{12}(\cc^2)$ so that $p = f^2+g^3+h^7$. 

By taking $d = 2s$, $e_1 = 2$, $e_2 = \dots = e_{s-1} = 1$ and $m=0$,
in Theorem \ref{omnibus}, we obtain an alternative to the dangling
term  `` $\la x^{2s}$'' in \eqref{E:Sylvcaneven}.
\begin{corollary}\label{6-22}
A general binary form $p$ of even degree $2s$ can be written as
\begin{equation}\label{E:sylvalt}
p(x,y) = (\al_0 x^2 + \be_0 x y + \ga_0 y^2)^s  
+ \sum_{j=1}^{s-1} (\al_j x + \be_j y)^{2s}.
\end{equation}
\end{corollary}

A different generalization of Theorem \ref{Sylvcan} focuses on the
number of summands.

\begin{theorem}\label{Sylvgen}
A general binary form of degree $uv$ can be written as a sum of
$\lceil \frac{uv+1}{u+1} \rceil$ $v$-th powers of binary forms of
degree $u$.
\end{theorem}

Cayley proved that,
after an invertible linear change of variables $(x,y) \mapsto (X,Y)$,
a general binary 
quartic can be written as $X^4 + 6 \la X^2Y^2 + Y^4$. 
 There are two natural ways to generalize
this to higher even degree, and almost 100 years ago, Wakeford
\cite{W2,W} did both. 
\begin{theorem}[Wakeford's Theorem]\label{missing}
After an invertible linear change of variables, a general $p \in H_d(\cc^n)$ can
be written so that the coefficient of each $x_i^d$ is 1 and the
coefficient of each $x_i^{d-1}x_j$ is 0.
\end{theorem}
There are $N(n,d) - n^2$ unmentioned monomials above, and when
combined with the $n^2$ coefficients in the change of variables, the
constant count is correct for a canonical form.
Wakeford was also interested in knowing {\it which} sets of
$n(n-1)$ monomials can be eliminated by a change of variables, and we
are able to settle this for binary forms in Theorem \ref{quarticgen}.
(Theorem \ref{missing} was independently discovered by
Guazzone \cite{Guaz} in 1975, as an attempt to generalize the canonical
form $X^3 + Y^3 + Z^3 + 6\la XYZ$ for $H_3(\cc^3)$. Babbage \cite{B}
subsequently observed that this can be proved by the Lasker-Wakeford
Theorem, without noting that Wakeford had already done so in
\cite{W}.)

The second generalization of $X^4 + 6 \la X^2Y^2 + Y^4$ will not be
pursued here; see \cite[Corollary 4.11]{ER}. A canonical form
for binary forms of even degree $2s$ is given by 
\begin{equation}\label{E:sylwake}
\sum_{k=1}^s \ell_k^{2s}(x,y) + \la \prod_{k=1}^s \ell_k^2(x,y),
\quad \ell_k(x,y) = \al_k x + \be_k y.
\end{equation}
This construction is due to Sylvester \cite{S2} for $2s=4,8$. His
methods failed for $2s=6$, but Wakeford was able to prove it in
\cite{W2}. The full version of \eqref{E:sylwake} 
is proved in \cite[p.408]{W}, where  Wakeford
notes that ``the number of ways this reduction can be performed is
interesting'', citing ``3,8,5'' for $2s=4,6,8$.

The non-trivial study of canonical forms was initiated by  Clebsch's 1861
discovery (\cite{Cleb}, see e.g. \cite[pp.50-51]{Ge} and
\cite[pp.59-60]{R2})  that, despite the fact that
$N(3,4) =5 \times N(3,1)$, a 
general ternary quartic cannot be written as a sum of five fourth
powers of linear forms. This was early evidence that constant-counting can fail.
But $N(3,4)$ is also equal to $1 \times N(3,2)+ 3 \times N(3,1)$,
and ternary quartics {\it do} satisfy an alternative canonical form as a
mixed sum of powers. 
\begin{theorem}\label{notclebsch}
A general ternary quartic $p \in H_4(\cc^3)$ can be written as 
\begin{equation}\label{E:notclebsch}
p(x_1,x_2,x_3) =  q^2(x_1,x_2,x_3) + \sum_{k=1}^3 \ell_k^4(x_1,x_2,x_3),
\end{equation}
where $q \in H_2(\cc^3)$ and $\ell_k \in H_1(\cc^3)$.
\end{theorem}

As an alternative generalization of canonical forms, one might 
also consider polynomial maps $F: S \mapsto H_d(\cc^n)$, where $S$ is an 
$N$-dimensional subspace of some $\cc^M$. 
In the simplest case, for binary quadratic forms, observe that
the coefficient of $x^2$ in
\begin{equation}\label{E:3}
(t_1 x + t_2 y)^2 + (it_1 x+ t_3 y)^2,
\end{equation}
is 0, so \eqref{E:3} is not canonical. This is essentially the only
kind of exception.
\begin{theorem}\label{hyperplane}
Suppose $(c_1,c_2,c_3,c_4) \in \cc^4$, and it is not true that $c_3 =
\ep c_1$ and $c_4 = \ep c_2$ for $\ep \in \{\pm i\}$.  Then for
general $p \in H_2(\cc^2)$, there exists $(t_1,t_2,t_3,t_4) \in \cc^4$
satisfying $\sum_{j=1}^4 c_jt_j = 0$ and such that
\begin{equation}\label{E:4}
p(x,y) = (t_1 x + t_2 y)^2 + (t_3 x+ t_4 y)^2.
\end{equation}
In the exceptional case,  there exists $(x_0,y_0)$ so that for all
feasible choices of $t_j$, $p(x_0,y_0) = 0$.
\end{theorem}
Another alternative version of \eqref{E:Sylvcaneven} is the following
conjecture, which can be verified up to degree 8.
\begin{conjecture}\label{zerosum}
A general binary form $p$ of even degree $2s$ can be written as
\begin{equation}\label{E:sylvalt2}
p(x,y) = \sum_{j=1}^{s+1} (\al_j x + \be_j y)^{2s}, \quad \text{where}\ 
\sum_{j=1}^{s+1} (\al_j + \be_j) = 0.
\end{equation}
\end{conjecture}

\subsection{Outline}

Here is an outline of the paper. In Section 2, we introduce notation
and definitions. The definition of canonical form is the classical
one and roughly parallels that in Ehrenborg-Rota \cite{ER}, an
important updating of this subject about 20 years ago. Our point of
view is considerably more elementary in many respects than \cite{ER},
but uses the traditional criterion: 
A polynomial map $F: \cc^N \mapsto H_d(\cc^n)$ is a {\it canonical form} if
a general $p \in H_d(\cc^n)$ is in the range; this occurs if and
only if there is at least one point $u \in \cc^N$ so that $\{\frac{\partial
  F}{\partial t_j}(u)\}$   spans $H_d(\cc^n)$. (See Corollary
\ref{jake}.)
This leads to 
immediate non-constructive proofs of Theorems \ref{uppertri},
\ref{sextican}, \ref{missing} and \ref{notclebsch}, and a somewhat
more complicated proof of Theorem \ref{quarticgen}, which answers
Wakeford's question about missing monomials for binary forms.

In Section 3, we discuss classical apolarity and its implications
for canonical forms. (Apolarity methods become more complicated when
a component of a canonical form comes from a restricted set of
monomials.) A generalization 
of the classical Fundamental   Theorem of Apolarity from \cite{R3}
allows us to identify a class of bases for $H_d(\cc^n)$ which give a
non-constructive proof of Theorem \ref{omnibus}, and 
hence Theorem \ref{Sylvcan}. A similar argument yields the proof of
Theorem \ref{Sylvgen}.  We also present Sylvester's
Algorithm, Theorem \ref{Sylvalg}, allowing for a constructive proof of
Theorem \ref{Sylvcan}.  We conclude
with a brief summary of  connections with the theorems of
Alexander-Hirschowitz and recent work on the rank of forms.

In Section 4 we discuss some special cases of Theorem
\ref{omnibus}. Sylvester's Algorithm is used in constructive
proof of Theorem \ref{omnibus} when $e_k \equiv 1$, in which case the
representation is unique. 
 We  give some other constructive proofs for $d \le 4$, and present
 numerical evidence regarding the number of representations 
in Corollary \ref{6-22} and a few other cases.  Using elementary
number theory, we show that, for each $r$, there are only finitely
many canonical forms \eqref{E:reveal} with $m=0$, and, up to degree
$N$, there are $N + \mathcal O(N^{1/2})$ such canonical forms in which
the $e_k$'s are equal. 

Section 5 discusses some familiar results on sums of two squares of
binary forms and canonical representations of quadratic forms as a sum
of squares of linear forms. This includes a
constructive proof of Theorem \ref{uppertri}, which
provides the groundwork for the proof of Theorem
\ref{slinkycan}.  We also
give a short proof of a canonical form which illustrates the classical 
result that a general ternary quartic is the sum of three squares of
quadratic forms.

In Section 6, we turn to forms in more than two variables and low degree,
give constructive proofs of 
Theorems \ref{reichcan} and \ref{slinkycan}, as well as the
non-canonical Theorem \ref{slowpoke}, which shows that {\it every} cubic in
$H_3(\cc^n)$ is a sum of at most $\frac{n(n+1)}2$ cubes of linear
forms.  Theorem \ref{reichcan} can be ``lifted'' to an ungainly
canonical form for quartics as a sum 
of fourth powers (see Corollary \ref{reichugly}), but not further to
quintics. Number theoretic considerations rule out a Reichstein-type
canonical form for quartics in 12 variables; see Theorem \ref{noreich}
for other instances of this phenomenon.

In Section 7, we offer a preliminary discussion of canonical forms
in which the domain of a polynomial map $F:\cc^M \mapsto H_d(\cc^n)$
is restricted to an 
$N$-dimensional subspace of $\cc^M$, of which Theorem \ref{hyperplane}
and  Conjecture \ref{zerosum}  are examples. 

The greatest debt of the author is due Richard Ehrenborg and
Gian-Carlo Rota for writing \cite{ER}.
Thanks to Dave Anderson and Julianna Tymoczko for organizing the 
Special Session on Geometric Commutative Algebra and Applications at
the March 2011 AMS Sectional Meeting in Iowa City, 
to Lek-Heng Lim for organizing the Minisymposium on Tensor Rank at 
the October 2011 SIAM Conference on Applied Algebraic Geometry in
Raleigh and to Eugene Mukhin for organizing ALGECOM5 in
Indianapolis in October 2011.  Invitations to speak at these
conferences provided an opportunity to present preliminary versions of this
material.

 I am indebted to Tony V\'arilly-Alvarado for extremely helpful
correspondence.  I also want to thank T. Y. Lam, Hal Schenck, Bernd
Sturmfels and Doron 
Zeilberger for their assistance. Special thanks go to the Center for
Advanced Study at 
UIUC, where the author was an Associate in the Fall 2011 semester, and
thereby free of teaching responsibilities. Finally, I want to thank
the referee for a careful reading of the manuscript and for making
many useful suggestions.

\section{Basic definitions, and proofs of Theorems 1.2, 1.5, 1.9 and
  1.10} 

Let $\mathcal I(n,d)$ denote the index set of monomials in $H_d(\cc^n)$:
\begin{equation}\label{E:fund}
\mathcal I(n,d) = \bigl\{(i_1,\dots,i_n): 0 \le i_k \in \mathbb Z,
  \quad  \sum_k i_k = d \bigr\}.
\end{equation}
Let $x^i = x_1^{i_1}\cdots x_n^{i_n}$ and $c(i) = \frac {d!}{\prod
  i_k!}$ denote the multinomial coefficient. If $p \in H_d(\cc^n)$,
then we  write
\begin{equation}\label{E:canon}
p(x_1,\dots,x_n) = \sum_{i \in \mathcal I(n,d)} c(i)a(p;i)x^i, \quad
a(p;i) \in \cc.
\end{equation}
We say that two forms are {\it distinct} if they are non-proportional,
and a set of forms is {\it honest} if the forms are pairwise distinct.
For later reference, recall Biermann's Theorem; see \cite[p.31]{R2}.
\begin{theorem}[Biermann's Theorem]\label{Biermann}
If $p \in H_d(\mathbb C^n)$ and $p \neq 0$, then there exists $i \in \mathcal
I(n,d)$ so that $p(i) \neq 0$. 
\end{theorem}

The easy verification of whether a formula is a canonical form for
$H_d(\cc^n)$ relies on a crucial alternative. A self-contained 
accessible proof is in \cite[Theorem 2.4]{ER}, for which
Ehrenborg and Rota thank M. Artin and A. Mattuck. For further
discussion of the underlying algebraic geometry,
see Section 9.5 in Cox, Little and O'Shea \cite{CLO}.
\begin{theorem}\label{alt} 
Suppose $M \ge N$ and $F: \cc^M \to \cc^N$ is a polynomial map; that is,
$$
F(t_1,\dots,t_M) = (f_1(t_1,\dots,t_M),\dots,f_N(t_1,\dots,t_M))
$$
where each $f_j \in \cc[t_1,\dots,t_M]$. Then either (i) or (ii) holds:

(i) The $N$ polynomials $\{f_j: 1 \le j \le N\}$ are algebraically
  dependent and $F(\cc^M)$ lies in some non-trivial variety $\{P =
  0\}$ in $\cc^N$. 
 
(ii) The $N$ polynomials $\{f_j: 1 \le j \le N\}$ are algebraically
  independent and $F(\cc^M)$ is dense in $\cc^N$.

The second case occurs if and only there is a point $u
\in \cc^M$ at which the Jacobian matrix
$\bigl[ \frac {\partial f_i}{\partial t_j}(u) \bigr]$
has full rank. 
\end{theorem}

When $M=N = N(n,d)$, we may interpret such an  $F$ as a map
from $\cc^N$ to $H_d(\cc^n)$ by indexing $\mathcal I(n,d)$ as $\{i(k): 
1\le k \le N\}$  and
making the interpretation in an abuse of notation that
\begin{equation}\label{E:abuse}
F(t;x) = \sum_{k=1}^N c(i(k))f_k(t_1\dots,t_N)x^{i(k)}.
\end{equation}

\begin{definition} A {\it canonical form for $H_d(\cc^n)$} is any 
polynomial map  $F : \cc^{N(n,d)} \mapsto H_d(\cc^n)$ in which $F$
satisfies  Theorem \ref{alt}(ii). 
\end{definition}
That is, $F$ is a canonical form if and only if $N = N(n,d)$ and for
a general $p \in   H_d(\cc^n)$, there exists $t \in \cc^N$ so that
$p(x) = F(t;x)$. The significance of this choice of $N$ is that it is
the smallest possible value.
In the rare cases where $F$ is surjective, we say that the
canonical form is {\it universal}.

By translating the definitions and using \eqref{E:fund} and
\eqref{E:abuse}, we obtain an immediate corollary of Theorem \ref{alt}:

\begin{corollary}\label{jake}
The polynomial map  $F : \cc^N \mapsto H_d(\cc^n)$ is a canonical form
if and only if there exists  $u \in \cc^n$ so that
$\{\frac{\partial F}{\partial t_j}(u)\}$  spans $H_d(\cc^n)$.
\end{corollary}

We shall let  $J := J(F;u)$ denote the span of the forms $\{\frac{\partial
  F}{\partial t_j}(u)\}$. In any particular case, the determination of
whether $J = H_d(\cc^n)$ amounts to the computation of the determinant
of an $N(n,d)\times N(n,d)$ matrix. As much as possible in this paper,
we give proofs which can be checked by hand, by making a judicious
choice of $u$ and ordering of the monomials in  $H_d(\cc^n)$, showing
sequentially that they all lie in $J$.

Classically, the use of the term  ``canonical form'' has been limited
to cases in 
which $F(t;x)$ has a natural interpretation as a combination of  forms in
$H_d(\cc^n)$, such as a sum of powers of linear forms, or as a result
of a linear change of variables. It seems odd
that canonical forms are perceived as rare, since a ``general''
polynomial map from $\cc^N \mapsto H_d(\cc^n)$ is a canonical
form. (This is an observation which
goes back at least to \cite{Rich}.) For example, 
if $\{f_j(x)\}$ is a basis for $H_d(\mathbb C^n)$, then
\begin{equation}\label{E:basis2}
F(t;x) = \sum_{j=1}^{N} t_j f_j(x)
\end{equation}
should be (but usually isn't) considered a canonical form. In
particular, \eqref{E:canon} with $f_j(x) = c(i_j)x^{i_j}$ is itself a
canonical form.

The following computation will occur repeatedly. If $es=d$, then
\begin{equation}\label{E:trivial}
g = \sum_{i_j \in \mathcal I(n,e)}t_j x^{i_j} \implies \frac{\partial
  g^s}{\partial t_j} = s x^{i_j}g^{s-1}.
\end{equation}
If $g$ is specialized to be a monomial, then all these partials will also
be monomials.

\begin{proof}[Non-constructive proof of Theorem \ref{uppertri}]
Given \eqref{E:upper}, let 
\begin{equation*}
\ell_k(x) =  \sum_{m=k}^n t_{k,m}x_m, \quad F(x) = \sum_{k=1}^n \ell_k^2(x).
\end{equation*}  Then
$\frac{\partial F}{\partial t_{k,m}} = 2x_m\ell_k$.  Set $t_{k,m}
= \de_{k,m}$, so that $\ell_k = x_k$ and 
$\frac{\partial F}{\partial t_{k,m}}  = 2x_kx_m$. Since $1
\le k \le m \le n$, all monomials from $H_2(\cc^n)$ appear in $J$.
\end{proof}
\begin{proof}[Non-constructive proof of Theorem  \ref{sextican}]
Suppose
\begin{equation}\label{E:nitty}
\begin{gathered}
p(x,y) = f^2(x,y) + g^3(x,y): \\
f(x,y) = t_1 x^3 + t_2 x^2y + t_3 xy^2 + t_4 y^3, \quad
g(x,y) = t_5 x^2 + t_6 xy + t_7 y^2.
\end{gathered}
\end{equation}
Then by \eqref{E:trivial}, the partials with respect to the $t_j$'s are:
\begin{equation*}
2x^3f,\ 2x^2yf,\  2xy^2f,\ 2y^3f; \quad 3x^2g^2,\ 3xyg^2,\ 3y^2g^2 .
\end{equation*}
Upon specializing at $f = x^3, g= y^2$, these become: 
\begin{equation*}
2x^6,\ 2x^5y,\  2x^4y^2,\ 2x^3y^3; \quad 3x^2y^4,\ 3xy^5\ , 3y^6.
\end{equation*}
It is then evident that  $J = H_6(\cc^2)$.
\end{proof}

\begin{proof}[Non-constructive proof of Theorem \ref{missing}]
Let $\mathcal L \subset \mathcal I(n,d)$ consist of all $n$-tuples
except the permutations of $(d,0,\dots,0)$ and $(d-1,1,\dots,0)$ and
let $X_i = \sum_{j=1}^n \al_{ij}x_j$. The assertion is that, with the
$(N(n,d) - n - \binom n2) + n^2 = N(n,d)$ parameters $t_{\ell}$ and $\al_{ij}$, 
\begin{equation}\label{E:wak}
\sum_{i=1}^n X_i^d + \sum_{\ell \in \mathcal L}
t_{\ell}X_1^{\ell_1}\cdots X_n^{\ell_n}.
\end{equation}
is a canonical form. Evaluate the partials at the point where $X_i
= x_i$ and $t_{\ell} = 0$:  they are $dx_jx_i^{d-1}$ (for $\al_{ij}$) 
and $x^{\ell}$ (for $t_{\ell}$). Taking $1 \le i, j \le n$ and $\ell \in
\mathcal L$, we see that $J$ contains all monomials in $H_d(\cc^n)$.
\end{proof}
As a special case (used later in Theorem \ref{6}), we obtain the
familiar result that 
after appropriate linear changes of variable, 
a general binary quartic may be written as $x^4 + 6\la x^2y^2 +
y^4$. It is classically known (see \cite[\S211]{EQ}) the choice of
$\la$ is not unique: in fact, after appropriate linear changes of
variable, $x^4 + 6\la x^2y^2 + y^4$ can be written as  $x^4 + 6\mu x^2y^2 + y^4$ 
for $\mu \in \{\pm \la, \pm \frac{1-\la}{1+3\la},
\pm \frac{1+\la}{1 - 3\la}\}$. 

Wakeford asserts that Theorem \ref{missing} is also true with
$x_i^{d-1}x_j$ replaced by $x_i^{d-r}x_j^r$ (evidently when $r \neq
\frac d2$), but his proof seems sketchy. He also gives necessary
conditions for sets of $n(n-1)$ monomials which may be omitted, and
these are hard to follow as well. Below, we answer his question in the
binary case: in the only  two excluded cases below,
\eqref{E:quarticgen} has a square factor, and so cannot be canonical.
\begin{theorem}\label{quarticgen}
Let $\mathcal B = (m_1,m_2,n_1,n_2)$ be four distinct
integers in $\{0,\dots,d\}$ so that 
$\{m_1,m_2\} \neq \{0,1\},\{d-1,d\}$. Then, after 
an invertible linear change of variable, a
general binary form $p$ of degree $d$ can be written as 
\begin{equation}\label{E:quarticgen}
p(x,y) = x^{d-n_1}y^{n_1} +x^{d-n_2}y^{n_2}  + \sum_{k\notin \mathcal B}
t_k x^{d-k}y^k
\end{equation}
for some $\{t_k\} \subset \cc$. 
\end{theorem}
\begin{proof}
Writing $(x,y) \mapsto (\al_1 x + \al_2 y, \al_3 x + \al_4 y):=
(X,Y)$, we have
\begin{equation}\label{E:wakegen}
F =X^{d-n_1}Y^{n_1} +X^{d-n_2}Y^{n_2}  + \sum_{k\notin \mathcal B}  t_k X^{d-k}Y^k. 
\end{equation}
Evaluate the partials of \eqref{E:wakegen} at
$(\al_1,\al_2,\al_3,\al_4) = (1,0,0,1)$ (so 
$X=x, Y=y$) and $t_k = 1$ (note the difference with the previous
proof, in which $t_k = 0$). The $d-3$ partials with respect to the $t_k$'s are 
simply $x^{d-k}y^k$, $k\notin \mathcal B$, so these are in $J$.  Further,
\begin{equation}\label{E:14}
\begin{gathered}
\frac{\partial F}{\partial \al_1} = \sum_{j \neq  m_1,m_2}(d-j)
x^{d-j}y^j, \quad \frac{\partial F}{\partial \al_4} = \sum_{j \neq  m_1,m_2}j
x^{d-j}y^j.
\end{gathered}
\end{equation}
Since most monomials used in \eqref{E:14} are already in $J$, it follows
that $J$ also contains  
\begin{equation}\label{E:142}
(d-n_1)x^{d-n_1}y^{n_1} + (d-n_2)x^{d-n_2}y^{n_2},\quad n_1x^{d-n_1}y^{n_1}
+ n_2x^{d-n_2}y^{n_2},
\end{equation}
and since $(d-n_1)n_2 \neq
(d-n_2)n_1$, \eqref{E:142} implies that $x^{d-n_j}y^{n_j} \in J$ for
$j=1,2$. To this point, we have shown that $J$ contains all
monomials from $H_d(\cc^2)$ except for $x^{d-m_j}y^{m_j}$, where $m_1 < m_2$. 
The two remaining partial derivatives are 
\begin{equation}\label{E:par23}
\begin{gathered}
\frac{\partial F}{\partial \al_2} = 
\sum_{j \neq  m_1,m_2}(d-j) x^{d-j-1}y^{j+1},\quad
\frac{\partial F}{\partial \al_3} =  
\sum_{j \neq  m_1,m_2} j x^{d-j+1}y^{j-1}, 
\end{gathered}
\end{equation}
and so $J$ contains as well the forms in \eqref{E:par23} of the shape
$c_1 x^{d-m_1}y^{m_1}+ c_2  x^{d-m_2}y^{m_2}$. We need to distinguish
a number of cases. If $m_1 = 0, m_2 = d$, then these forms are
$y^d, x^d$. If $m_1 = 0$ and $2 \le m_2
\le d-1$, then these forms are $(d-m_2) x^{d-m_2}y^{m_2}$ and $x^d +
(m_2+1)x^{d-m_2}y^{m_2}$, and similarly when $1 \le m_1 \le d-2$ and
$m_2=d$. (Recall that we have excluded the cases $(m_1,m_2) = (0,1)
and (d-1,d)$. In the remaining cases,  
 $1 \le m_1 < m_2 \le d-1$. If $m_2 = m_1 + 1$, then these forms are 
$(d - (m_1-1))x^{d-m_1}y^{m_1}$ and $(m_2+1)x^{d-m_2}y^{m_2}$.  Finally, if  
$m_2 > m_1 + 1$, then all four terms appear, and the forms are
\begin{equation}\label{E:par232}
\begin{gathered}
(d-m_1+1)x^{d-m_1}y^{m_1} + (d-m_2+1)x^{d-m_2}y^{m_2},\\
(m_1+1)x^{d-m_1}y^{m_1} + (m_2+1)x^{d-m_2}y^{m_2}. 
\end{gathered}
\end{equation}
In each of the cases, linear combinations of the forms produce the
missing monomials, so  $J = H_d(\cc^2)$. 
\end{proof}
\begin{remark}
By writing $p(x,y) = \prod_k(x +\al_k
y)$, it follows  from Theorem \ref{missing} that, for a general set of
$d$ complex numbers $\al_k$, there exists a M\"obius transformation
$T$ so that  
\begin{equation}
\sum_{k=1}^d T(\al_k) = 0,\quad  \sum_{k=1}^d T(\tfrac 1{\al_k}) = 0,
\quad \prod_{k=1}^d T(\al_k) = 1. 
\end{equation}
\end{remark}

\begin{proof}[Non-constructive proof of Theorem  \ref{notclebsch}]
Write \eqref{E:notclebsch} as $F(x;t)$, where
\begin{equation*}
\begin{gathered}
q(x_1,x_2,x_3) = t_1x_1^2 + t_2 x_2^2 + t_3 x_3^2 + t_4 x_1x_2 + t_5
x_1x_3 + t_6 x_2x_3, 
\\
\ell_k(x_1,x_2,x_3) = t_{k1}x_1 + t_{k2}x_2 +t_{k3}x_3. 
\end{gathered}
\end{equation*}
Evaluate the partials at: 
$q =  x_1x_2 +  x_1x_3 + x_2x_3$ and $(\ell_1,\ell_2,\ell_3) =
(x_1,x_2,x_3)$.
Then $\frac{\partial F}{\partial t_{k\ell}} = 4x_{\ell}x_k^3$, so
$x_i^4, x_i^3x_j \in J$;  since  $\frac{\partial F}{\partial
  t_1}=2x_1^2q = 2x_1^2(x_1x_2 +  x_1x_3 + x_2x_3)$, it follows that
$x_1^2x_2x_3 \in J$, similarly, by considering $\frac{\partial F}{\partial
  t_2}$ and $\frac{\partial F}{\partial  t_3}$, it follows that
$x_1x_2^2x_3,  x_1x_2x_3^2$ are in $J$. Finally,
$\frac{\partial F}{\partial t_4} = 2x_1x_2q =2 x_1x_2(x_1x_2 +  x_1x_3
+ x_2x_3)$, and so now $x_1^2x_2^2 \in J$. Similarly,  by considering
$\frac{\partial F}{\partial  t_5}$ and $\frac{\partial F}{\partial
  t_6}$, it follows that $x_1^2x_3^2,
x_2^2x_3^2$ are also in  $J$, and this accounts for all monomials in
$H_4(\cc^3)$. 
\end{proof}

Other applications of Corollary \ref{jake}  to canonical forms
can be found in \cite{ER}, including interpretations of the older
results in \cite{Rich} and \cite[pp.265-269]{T}.  

\section{Apolarity and proofs of Theorems 1.1, 1.6 and 1.8}

Using the notation of \eqref{E:fund} and \eqref{E:canon}, for $p, q
\in H_d(\cc^n)$, define  the following bilinear form: 
\begin{equation}\label{E:ip}
[p,q] = \sum_{i \in \mathcal I(n,d)} c(i)a(p;i)a(q;i).
\end{equation}

Recall two basic notations.
For $\al \in \cc^n$, define $(\al\cdot)^d \in H_d(\cc^n)$ by 
\begin{equation}
(\al\cdot)^d(x) = (\al\cdot x)^d = \biggl(\sum_{j=1}^n\al_j x_j\biggr)^d =
  \sum_{i \in \mathcal I(n,d)} c(i) \al^i x^i.
\end{equation}
Define the
differential operator $f(D)$ for $f \in H_e(\cc^n)$ in the usual way by
\begin{equation}\label{E:diff}
f(D) = \sum_{i \in \mathcal I(n,e)} c(i)a(f;i)\left(
  \tfrac{\partial}{\partial x_1}\right)^{i_1}\cdots
\left(\tfrac{\partial}{\partial 
    x_n}\right)^{i_n}.  
\end{equation}
It follows immediately that for $\al \in \cc^n$,
\begin{equation}\label{E:dual}
[p,(\al\cdot)^d] =  \sum_{i \in \mathcal I(n,d)} c(i)a(p;i)\al^i = p(\al).
\end{equation}
If $i\neq j \in \mathcal I(n,d)$, then $i_k > j_k$ for some $k$, so $D^{i}x^{j} =
0$; otherwise $D^ix^i = \prod_k (i_k)! = 
d!/c(i)$. Suppose $p,q \in H_d(\cc^n)$.
 Bilinearity and \eqref{E:diff}  imply the classical result that
\begin{equation}\label{E:difflin}
\begin{gathered}
p(D)q = \sum_{i \in \mathcal I(n,d)} c(i)a(p;i)D^i \biggl(\sum_{j \in
  \mathcal I(n,d)} c(j)a(q;j)x^j\biggr) =
\\ \sum_{i \in \mathcal I(n,d)}  \sum_{j \in \mathcal I(n,d)}
c(i)c(j)a(p;i)a(q;j)D^ix^j 
 = \sum_{i \in \mathcal I(n,d)}  c(i)c(i)a(p;i)a(q;i)D^ix^i 
\\ = \sum_{i \in \mathcal I(n,d)}  c(i)^2a(p;i)a(q;i) \frac {d!}{c(i)} =
d![p,q] = d![q,p] = q(D)p.
\end{gathered}
\end{equation}

\begin{definition}
If $p \in H_d(\cc^n)$ and $q \in H_e(\cc^n)$, then $p$ and $q$ are
{\it apolar} if $p(D)q = q(D)p = 0$. 
\end{definition}

Note that if $d=e$, then $p$ and $q$ are apolar if and only if
$[p,q]=0$ and if $d > e$, say, then the equation $p(D)q = 0$ is
automatic, so only $q(D)p = 0$ need be checked. 
By \eqref{E:dual}, $p$ is apolar to $(\al\cdot)^d$ if and only
if $p(\al) = 0$.  

The following lemma is both essential and trivial.
\begin{lemma}\label{dual}
Suppose $X = span(\{h_j\}) \subseteq H_{d}(\cc^n)$. 
Then $X =  H_{d}(\cc^n)$ if and only if there is no  $0 \neq p \in H_d(\cc^n)$
which is apolar to each of the $h_j$'s.
\end{lemma}
From this point of view, Theorem \ref{LW} is a direct consequence
of Corollary \ref{jake}:
\begin{theorem}[Lasker-Wakeford]\label{LW}
If $F: \cc^N \to H_d(\cc^n)$, then $F$ is a canonical form if and only 
if there is a point $u$ so that there is no non-zero form $q \in
H_d(\cc^n)$ which is
apolar to all $N$ forms $\{\frac{\partial F}{\partial t_k}(u)\}$. 
\end{theorem}

The attribution ``Lasker-Wakeford" (for \cite{Lasker,W}) is taken from
\cite{T}:  H. W. Turnbull 
(1885-1961) was one of the last
practicing invariant theorists who had been trained in the pre-Hilbert
approach, see \cite[pp.231-232]{Fi}. (His text \cite{T} is  a 
Rosetta Stone for understanding the 19th century approach to algebra
in more modern terminology.)  Turnbull referred to
Theorem \ref{LW} as ``paradoxical and very curious''.  
E. Lasker (1868-1941) received his Ph.D. under M. Noether at
G\"ottingen in 1902. He is probably better
known for being the world chess champion for 27 years (1894-1921),
spanning the life of E. K. Wakeford (1894-1916). 
J. H. Grace,  Wakeford's professor at Oxford, edited the second half
of his thesis into the article \cite{W} and also wrote a memorial
article \cite{Gr} for him in 1918:
\begin{quote}
``He [EKW] was slightly wounded early in 1916, and soon after coming home
was busy again with Canonical Forms.... [H]e discovered a paper of
Hilbert's which contained the very theorem he had long been in want of
-- first vaguely, and later quite definitely. This was in March; April
found him, full of the most joyous and reverential admiration for the
great German master, working away in fearful haste to finish the
dissertation ...  He returned to the front in June and was killed in
July.... He only needed a chance, and he never got it.'' 
\end{quote}

The following properties are easily established; see, e.g.,
\cite{R2,R3} for proofs.

\begin{theorem}\label{ip}
\

\noindent (i)   If $e \le d$ and $f \in H_e(\cc^n)$, $g \in H_{d-e}(\cc^n)$ and $p
\in H_d(\cc^n)$, then 
\begin{equation}\label{E:facto}
d![fg,p] = (fg)(D)p = f(D)g(D)p = e![f,g(D)p].
\end{equation}
Thus, $p$ is
apolar to every multiple of $g$ in $H_d(\cc^n)$ if and only if $p$ and
$g$ are apolar.

\noindent (ii)  If $p
\in H_d(\cc^n)$, then  $\tfrac 1d \tfrac{\partial p}{\partial x_j}(\al) =
[p,x_j(\al\cdot)^{d-1}]$. Thus, 
$p$ is apolar to $(\al\cdot)^{d-1}$ 
if and only if $p$ is singular at $\al$. More generally,  $p$ is
apolar to $(\al\cdot)^{d-e}$ if and only if $p$ vanishes to $e$-th
order at $\al$. 

\noindent (iii)  If $ e \le d$ and  $g \in H_{d-e}(\cc^n)$, then 
$g(D)(\al\cdot)^d = \frac  {d!}{e!} g(\al)(\al\cdot)^e$.

\end{theorem}

 Suppose $F(t;x)$ contains $h^s$ as a summand, where $h(x) = \sum_{\ell
  \in \mathcal I(n,e)}t_{\ell}x^{\ell}$, and suppose that no
$t_{\ell}$ occurs elsewhere in $F(t;x)$. If $p$ is apolar to each
partial of $F$, then it will be apolar to $\frac{\partial
  F}{\partial t_{\ell}} = sx^{\ell}h^{s-1}$ by  \eqref{E:trivial}.
 Since this is true for
every $\ell \in \mathcal I(n,e)$, it follows from (i) that $p$ is
apolar to $h^{s-1}$. It is critical to note that this observation
requires that each of the monomials of degree $e$ appear in $h$, and
does not apply if $h$ is defined as a sum from a restricted set of
monomials.
 
We are now able to give a short proof of the 
``Second main theorem on apolarity'' from \cite{ER}, which was
not concerned with preserving the constant-count. 
\begin{theorem}\label{ER}
Suppose $j_{\ell} = (j_{\ell,1}, \dots, j_{\ell,m}), 1 \le \ell \le r$,
are $m$-tuples of non-negative integers, and suppose positive integers $d_k$,
$1 \le k \le m$, and $d$ are chosen so that
\begin{equation}\label{E:homre}
u_{\ell}:= d - \sum_{k=1}^m j_{\ell,k}d_k \ge 0
\end{equation} 
for each $\ell$. Fix forms $q_{\ell} \in H_{u_\ell}(\cc^n)$ and for $f_k \in
H_{d_k}(\cc^n)$, define
\begin{equation}\label{E:ER}
F(f_1,\dots,f_m) = \sum_{\ell=1}^r q_{\ell}(x)f_1^{j_{\ell,1}}\cdots  f_m^{j_{\ell,m}}.
\end{equation}
Let $F_j:= \frac {\partial F}{\partial f_j}$. Then
a general $p \in H_d(\cc^n)$ can be written as \eqref{E:ER} if and
only if there exists a specific  $\bar f = (\bar f_k)$ so 
that no non-zero $p \in H_d(\cc^n)$ is apolar to each $F_j(\bar f)$,
$1 \le j \le m$. 
If, in addition, 
\begin{equation}\label{E:addem}
\sum_{k=1}^m N(n,d_k) = N(n,d),
\end{equation}
then \eqref{E:ER} is a canonical form.
\end{theorem}

\begin{proof}
Let
\begin{equation}\label{E:39}
 f_j(x) = \sum_{i_v \in \mathcal I(n,d_j)} t_{j,v} x^{i_v}. 
\end{equation}
By Theorem  \ref{alt}, \eqref{E:homre} and Lemma \ref{dual}, 
\eqref{E:ER} represents
general $p \in  H_d(\cc^n)$ if and only if there is some $\bar f$ so that
there is no non-zero form in $p \in H_d(\cc^n)$ which is apolar to
each $\frac {\partial F}{\partial t_{j,v}}(\bar f) = d_k
x^{i_v}F_j(\bar f)$, or by Theorem \ref{ip}(i), to each $F_j(\bar f)$. 
 The constant count is checked by \eqref{E:addem}.
\end{proof}

By Theorem \ref{ip}(ii) and Theorem \ref{ER}, 
 $F = \sum_{k=1}^r (\al_k\cdot)^d$ is a canonical form if and only if
there exist $r$ points $\bar\al_k \in \cc^n$ at which no
non-zero form  $p \in H_d(\cc^n)$ is singular.  This result is
classical, and goes back to Clebsch \cite{Cleb}; see 
also \cite[Theorem 4.2]{ER}.  A particularly deep result of
Alexander and Hirschowitz \cite{AH} from the early 
1990s states that a general form in $H_d(\cc^n)$, $d \ge 3$, may be written as a
sum of $\lceil \frac 1n N(n,d) \rceil$ $d$-th powers of linear forms,
except when  $(n,d) = (5,3), (3,4), (4,4), (5,4)$, when an extra
summand is needed.  (For much more on this, see \cite[Lecture 7]{Ge}, 
\cite[Corollary 1.62]{IK}, \cite[Chapter 15]{Lan} and
\cite[Theorem 0.2]{RanSch}; for a brief exposition of the proof, 
 see \cite[Chapter 15]{Lan}.) These references also discuss the exceptional
 examples, which  were all known in the 19th century. The expression
 of forms as a sum of powers of  forms is currently a very active
 area of interest; see the references above as well as \cite{CCG},
 \cite{FOS} and \cite{LT}.

The Fundamental Theorem of Apolarity (see \cite{R3} for a 
history) states that if $f$ is irreducible and $p \in H_d(\cc^n)$, 
then $f$ and $p$ are apolar  if and only if $p$  can be written as a sum of
terms of the form $(\al_j\cdot)^d$, where
 $f(\al_j) = 0$.  This was generalized in \cite{R3}.

\begin{theorem}\cite[Theorem 4.1]{R3}\label{adv}
Suppose $q \in H_e(\cc^n)$ factors as $\prod_{j=1}^r q_j^{m_j}$ into a product
of powers of distinct irreducible factors and suppose $p \in
H_d(\cc^n)$. Then $q(D)p = 0$ if and only if there exist $\al_{jk} \subset \{
q_j(\al) = 0\}$, and $\phi_{jk} \in H_{m_j-1}(\cc^n)$ such that 
\begin{equation*}
p = \sum_{j=1}^r\left(\sum_{k=1}^{n_j}
  \phi_{jk}(\al_{jk}\cdot)^{d-(m_j-1)}\right). 
\end{equation*}
\end{theorem}

The application of apolarity to binary forms is particularly simple,
because zeros correspond to factors. If $e=d+1$, then  $q(D)p =
0$ for every $p \in H_d(\cc^n)$, and we obtain the following result,
also found in \cite[Theorem 4.5]{ER}.
\begin{corollary}\label{cor}
Suppose $\{\al_j x + \be_j y: 1 \le j \le r\}$ is honest 
 and suppose $\sum_{j=1}^{r} m_j = d+1$. Then the following
set is a basis for $H_d(\cc^2)$:
\begin{equation}\label{E:mess}
\mathcal S = \left\{ x^k y^{m_j-1-k} (\be_j x - \al_j y)^{d-m_j+1}: 0
  \le k \le m_j-1,\quad 1  \le m_j \le r \right\}. 
\end{equation}
\end{corollary}
\begin{proof}
If $p$ is apolar to each term in \eqref{E:mess}, then $(\al_j x + \be_j
y)^{m_j}\ | \ p$ by Theorem \ref{ip}(ii). Thus $p=0$ by degree
considerations, and $\mathcal S$ has $d+1$ elements, so it is a basis.
\end{proof}
If each $m_j=1$, then Corollary \ref{cor} states that an honest set 
$\mathcal S= \{(\al_j x + \be_j y)^d\}$ of $d+1$ forms  is a basis for
$H_d(\cc^2)$. This is easily proved directly, since the representation
of $\mathcal S$ with respect to the basis  $\{\binom dj x^{d-j}y^j\}$,
$[\al_j^{d-k}\be_j^k]$, has Vandermonde determinant
\begin{equation}\label{E:vdm}
\prod_{1 \le i < j \le n} (\al_i\be_j - \al_j\be_i).
\end{equation}
Each product in \eqref{E:vdm} is non-zero because $\{(\al_j x + \be_j y)^d\}$
is honest. One implication of this independence is found in
\cite[Corollary 4.3]{R4}.
\begin{lemma}\label{lengthlemma}
If $p(x,y) \in H_d(\cc^2)$ has two honest representations
\begin{equation}\label{E:leng}
p(x,y) = \sum_{i=1}^m(\al_i x + \be_i y)^d = \sum_{j=1}^n(\ga_j x + \de_j y)^d
\end{equation}
and $m+n \le d+1$, then the representations are permutations of each other.
\end{lemma}
\begin{proof}
If \eqref{E:leng} holds, then $\{(\al_i x + \be_i y)^d, (\ga_j x +
\de_j y)^d\}$ is linearly dependent, which is impossible unless the
dependence is trivial. 
\end{proof}
 It follows immediately from Lemma \ref{lengthlemma} that the
representations \eqref{E:Sylvcanodd} and  \eqref{E:Sylvcaneven}, if they
exist for $p$, are unique. When $n \ge 3$, the linear dependence of
a set $\{(\al_j\cdot)^d\}$ depends on the geometry of the points as well
as the number (see the discussion of Serret's Theorem in \cite[p.29]{R2}.) 
Even for powers of binary forms of degree $e \ge 2$, there
are singular cases. It is not hard to show that a {\it general}
set of $(2k+1)$ $k$-th powers of quadratic forms is linearly
independent; however, 
for example, $(x^2-y^2)^2 + (2xy)^2 = (x^2+y^2)^2$. For much more on this,
see \cite{R7}.

\begin{proof}[Non-constructive proof of Theorem \ref{omnibus}]
For $1 \le k \le r$, write 
\begin{equation*}
f_k(x,y) = \sum_{\ell = 0}^{e_k} t_{k,\ell}x^{e_k - \ell} y^{\ell}.
\end{equation*}
By Corollary \ref{jake} and \eqref{E:trivial}, \eqref{E:reveal} is a
canonical form in the variables 
$\{t_j,t_{k,\ell}\}$ provided there is a point at which the partials 
\begin{equation*}
\{ \ell_j^d,  \ 1\le j \le m\} \cup \ \{
  x^{e_k-\ell}y^{\ell}f_k^{d/e_k-1}, \ 1
  \le \ell \le e_k,\quad 1 \le k \le r \} 
\end{equation*}
span $H_d(\cc^2)$. Let $f_k = \tilde \ell_k^{e_k}$, where
$\{\ell_1,\dots,\ell_m,\tilde \ell_1, \dots, \tilde \ell_r\}$ is
chosen to be honest. Then by \eqref{E:book}, the desired assertion follows
immediately from Corollary \ref{cor}. 
\end{proof}
\begin{proof}[Non-constructive proof of Theorem \ref{Sylvgen}]
Write $uv+1 = r(u+1)+s$. If $s=0$, then Theorem \ref{Sylvgen} is
simply a special case of Theorem \ref{omnibus} with $m=0$, $d=uv$ and
$e_k \equiv u$. Otherwise,  $1 \le s \le u$, so that $r+1 = \lceil
\frac{uv+1}{u+1} \rceil$. Let
\begin{equation*}
\begin{gathered}
F(\{\al_{ij}\}) = \sum_{i=1}^{r+1}f_i^v(x,y), \qquad
f_i(x,y) = \sum_{j=0}^u \al_{ij}x^{u-j}y^j.
\end{gathered} 
\end{equation*}
This is {\it not} a canonical form, as there are too many
constants. As before, $\frac{\partial F}{\partial \al_{ij}} =v x^{u-j}y^j 
f_i^{v-1}$.  We now specialize to  $f_i(x,y) = (ix - 
  y)^u$ and use the apolarity argument to show that $J =
  H_{uv}(\mathbb C^2)$.  Suppose $q \in  H_{uv}(\mathbb C^2)$ is
  apolar to each partial. Then by
Theorem \ref{ip}, it is apolar to $f_i^{v-1} = (ix -y)^{uv-u}$, and so $q$
vanishes to $u$-th order at $(i,-1)$ for $1 \le i \le r+1$. It follows that
$q$ is a multiple of $\prod_{i=1}^{r+1}(x + iy)^{u+1}$, and so $q=0$
by degree considerations. 

It is an exercise to show that $F$ can be converted to an
canonical form by requiring, say,  that $f_{r+1}$ only contain
monomials $x^{u-j}v^j$ for $0 \le j \le s-1$. 
\end{proof}

We present now Sylvester's Algorithm. For modern
discussions of this, along with Gundelfinger's generalization \cite{G},
which is not included here, see \cite[\S 5]{KR},  \cite{K1},\cite{K2},
\cite{K3}, \cite{R3} and \cite{R4}.
\begin{theorem}[Sylvester's Algorithm]\label{Sylvalg}
Let 
\begin{equation*}
p(x,y) = \sum_{j=0}^d \binom dj a_j x^{d-j}y^j
\end{equation*}
be a given binary form and suppose $\{\alpha_j x + \beta_j y\}$ is 
honest. Let 
\begin{equation*}
h(x,y) = \sum_{t=0}^r c_tx^{r-t}y^t = \prod_{j=1}^{r} (\beta_j x - \alpha_j y).
\end{equation*}
 Then there exist $\lambda_k\in \mathbb C$ so that  
\begin{equation*}
p(x,y) = \sum_{k=1}^r \lambda_k (\alpha_k x + \beta_k y)^d
\end{equation*}
if and only if
\begin{equation}\label{E:Syl}
\begin{pmatrix}
a_0 & a_1 & \cdots & a_r \\
a_1 & a_2 & \cdots & a_{r+1}\\
\vdots & \vdots & \ddots & \vdots \\
a_{d-r}& a_{d-r+1} & \cdots & a_d
\end{pmatrix}
\cdot
\begin{pmatrix}
c_0\\c_1\\ \vdots \\ c_r
\end{pmatrix}
=\begin{pmatrix}
0\\0\\ \vdots \\ 0
\end{pmatrix}.
\end{equation}
\end{theorem}
Theorem \ref{Sylvalg} can be put in the context of our previous discussion.
Let $A_r(p)$ denote the $(d-r+1) \times (r+1)$ Hankel matrix on the
left-hand side of 
\eqref{E:Syl}.
 If $h(D) = \prod_{j=1}^{r} (\beta_j
\frac{\partial}{\partial x} - \alpha_j \frac {\partial}{\partial y})$, then 
a direct computation shows that
\begin{equation}\label{E:syldifff}
h(D)p = \sum_{m=0}^{d-r} \frac{d!}{(d-r-m)!m!}\left( \sum_{i=0}^{d-r}
  a_{i+m}c_i \right) x^{d-r-m}y^m.
\end{equation}
It follows from \eqref{E:syldifff} that 
the coefficients of $h(D)p$ are thus, up to multiple, the rows of the
matrix product, so \eqref{E:Syl} is equivalent to $h(D)p=0$. In this way,
Theorem \ref{Sylvalg} follows from Theorem \ref{adv}.
Sylvester's algorithm can also
be visualized as seeking constant-coefficient linear recurrences
satisfied by $\{a_k\}$ and looking for the shortest one whose
characteristic equation has distinct roots; this is the proof given in
\cite{R4}. In this case, Gundelfinger's results handle the case when
the roots are not distinct.

\begin{proof}[Constructive proof of Theorem \ref{Sylvcan}]
Suppose $d = 2s-1$ is odd. The matrix $A_s(p)$ is $s \times (s+1)$ and
has a non-trivial null-vector. The corresponding $h$ (which can
be given in terms of the coefficients of $p$) has distinct factors unless its
discriminant vanishes. Thus for general $p \in H_{2s-1}(\cc^2)$,
Theorem \ref{Sylvalg} gives $p$ as a sum of $s$ $(2s-1)$-st powers of
linear forms.

If $d=2s$, the matrix $A_s(p)$ is square, and if $p$ is a sum of $s$
$2s$-th powers, then $\det A_s(p) = 0$. Conversely, if $\det A_s(p) =
0$ and the corresponding $h$ has distinct factors (which is generally
true), then $p$ is a sum of $s$ $2s$-th powers. If $M_1$
and $M_2$ are two square matrices and rank$(M_2) = k$, 
then $\det(M_1 + \la M_2)$ is a polynomial in $\la$ of degree $k$. In
particular,  if $q = 
(\al x + \be y)^{2s}$, then rank($H_s(q)) = 1$. Thus, in general,
there is a unique value of $\la$ and some matrix $M$ so that   $0 =
\det A_s(p-\la 
(\al x + \be y)^{2s}) = \det A_s(p) - \la \det M$. (When $\al x +
\be y = x$, $M$ is the  (1,1)-cofactor of $A_s(p)$.) In the special
case $\al x + \be y = x$, this proves Theorem \ref{Sylvcan}(ii).
The same argument
shows that for general $q \in H_{2s}(\cc^2)$, there exist $s+1$ values
of $\la$ so that $p - \la q$ is a sum of $s$ $2s$-th powers.
\end{proof}

In 1869, Sylvester \cite{S3} recalled his discovery of this algorithm
and its consequences.
\begin{quote}
``I  discovered and developed 
the whole theory of canonical binary forms for odd degrees, and, as far as yet 
made out, for even degrees too, at one evening sitting, with a decanter of port 
wine to sustain nature's flagging energies, in a back office in Lincoln's Inn 
Fields. The work was done, and well done, but at the usual cost of racking 
thought ---  a brain on fire, and feet feeling, or feelingless, as if
plunged in an  ice-pail. {\it That night we slept no more.}''
\end{quote}

\begin{example}\label{sylex}
This example of Sylvester's algorithm will be used in Example 4.1. Let
$p(x,y) =  2 x^3 + 3 x^2y - 21 x y^2 -41 y^3 =  
\binom 30 \cdot 2 \ x^3 + 
\binom 31 \cdot 1\ x^2y +  \binom 32 \cdot (-7)\ xy^2 +  \binom 33
\cdot(-41)\  y^3$
Since
\begin{equation*}
\begin{gathered}
\begin{pmatrix}
2 & 1 & -7\\
1 & -7 & -41\\
\end{pmatrix}
\cdot
\begin{pmatrix}
6\\-5\\1 
\end{pmatrix}
=
\begin{pmatrix}
0\\0\\0
\end{pmatrix},
\end{gathered}
\end{equation*}
we have $h(x,y) = 6x^2-5xy+y^2 = (2x - y)(3x-y)$. It now follows that
$p(x,y) = \la_1(x+2y)^3 + \la_2(x+3y)^3$, and a simple computation
shows that $\la_1 = 5, \la_2 = -3$.
\end{example}

 Lemma \ref{dual}, when applied to  Theorem
\ref{Biermann}, yields the following corollary. 
\begin{corollary}\label{altbier}
A basis for $H_d(\cc^n)$ is given by $\{(i \cdot )^d: i \in \mathcal I(n,d)\}$.
\end{corollary}
This in turn gives a very weak 
version of the Alexander-Hirschowitz Theorem,
\begin{corollary}\label{badah}
A general form in  $H_d(\cc^n)$
is a sum of $N(n,d-1) = \frac {nd}{n+d-1} \cdot \frac 1n N(n,d)$ $d$-th powers
of linear forms. 
\end{corollary}
\begin{proof}
Consider the sum 
\begin{equation*}
\sum_{\ell= 1}^{N(n,d-1)} (t_{\ell,1}x_1 + \cdots + t_{\ell,n}x_n)^d,
\end{equation*}
and apply Corollary \ref{jake} with $t_{\ell}$ specialized to
$i_{\ell} \in \mathcal I(n,d-1)$. Then $J$ contains
$x_k(i_{\ell}\cdot)^{d-1}$ for each 
$k,\ell$ and hence $x_kH_{d-1}(\cc^n) \subseteq J$ for each $k$, so $J =
H_d(\cc^n)$.  
\end{proof}

\section{Examples of binary canonical forms and the  proof
Theorem 1.7}

This section is devoted to special cases of Theorem \ref{omnibus}.
 First, in the special case $e_k=1$, we give a constructive proof
 showing uniqueness, 
 which gives a kind of interpolation between Sylvester's Theorem
  and  the representations of $H_d(\cc^2)$
 by \eqref{E:basis2} with a fixed  basis consisting of $d$-th powers,
 as in Corollary \ref{cor}.

\begin{corollary}\label{e=1}
Suppose $d \ge 1$, and $\{\ell_j(x,y) = \al_j x + \be_j y \}$ is a
fixed honest set of $m=d+1-2r$ 
linear forms. Then a general binary
$d$-ic form $p \in H_d(\cc^2)$ can be written uniquely as
 \begin{equation}\label{E:3.2}
p(x,y) = \sum_{j=1}^m t_j \ell_j(x,y)^d + \sum_{k=1}^r (t_{k1}x + t_{k2}y)^d.
\end{equation}
for suitable $t_{k1}, t_{k2} \in \cc$.
\end{corollary}
\begin{proof}
Let 
\begin{equation*}\label{E:killer}
f(x,y) = \prod_{j=1}^m(\be_j x - \al_j y).
\end{equation*}
Then $f(D)p$ has degree $d-m = 2r-1$ and by Theorem \ref{Sylvalg}
generally has a unique representation as a sum of $r$ $2r-1$-st powers
of linear forms, say
\begin{equation}\label{E:tg}
f(D)p = \sum_{k=1}^r (u_{k1} x + u_{k2} y)^{2r-1}.
\end{equation}
Further, it is generally true that $f(u_{k1},u_{k2}) \neq 0$. Let  
\begin{equation}\label{E:th}
q(x,y) = \frac{(2r-1)!}{d!} \sum_{k=1}^r 
\frac{(u_{k1} x + u_{k2} y)^d}{f(u_{k1},u_{k2})}.
\end{equation}
It follows from Theorem \ref{ip}(iii), \eqref{E:tg} and \eqref{E:th} that
$f(D)p = f(D)q$. Since $f$ has distinct factors, it then follows from
Theorem \ref{Sylvalg}  that there exist $t_j \in \cc$ so that 
\begin{equation*}
p(x,y) - q(x,y) = \sum_{j=1}^m t_j (\al_j x + \be_j y)^d.
\end{equation*}
Conversely, suppose $p$ has two different representations:
\begin{equation}\label{E:3.3}
 \sum_{j=1}^m t_j \ell_j^d(x,y) + \sum_{k=1}^r (t_{k1}x + t_{k2}y)^d =
 \sum_{j=1}^m \tilde t_j \ell_j^d(x,y) + \sum_{k=1}^r (\tilde t_{k1}x
 + \tilde t_{k2}y)^d.
\end{equation}
By combining the first sum on each side, \eqref{E:3.3} becomes a
linear dependence with $m+2r = d+1$ summands, which by Lemma \ref{lengthlemma}
must be trivial; thus, the representations in \eqref{E:3.3} are
essentially the same. 
\end{proof}
\begin{example}
Let $\ell_1(x,y) = x+y$ and $\ell_2(x,y) = -x + 3y$ and let
\begin{equation*}
 p(x,y) = -x^5 + 15 x^4y - 170x^3y^2 + 390 x^2y^3 - 505 x^2 y^3 + 483
y^5.
\end{equation*}
In an application of the last proof, $f(x,y)=(x-y)(3x+y)=
3x^2 -2xy-y^2$, and
\begin{equation*}
3 \frac {\partial^2 p }{\partial x^2}  - 2\frac {\partial^2 p}{\partial
  x\partial y} -  \frac {\partial^2 p}{\partial y^2} = 160 x^3 + 240
x^2 y - 1680 x y^2 - 3280 y^3.
\end{equation*}
Example \ref{sylex} implies that this expression equals
$400(x+2y)^3 - 240(x + 3y)^3$. Since $f(1,2) = -5$ and $f(1,3) = -12$,
it follows that
\begin{equation*}
\begin{gathered}
p(x,y) = \\  \frac{3!\cdot 400}{5!\cdot(-5)}(x + 2y)^5 +  \frac{3!\cdot
  (-240)}{5!\cdot(-12)}(x + 3y)^5 + t_1(x+y)^5 + t_2(-x + 3y)^5 = \\
 -4(x+2y)^5 +(x+3y)^5 + t_1(x+y)^5 + t_2(-x + 3y)^5
\end{gathered}
\end{equation*}
and it can be readily be computed that $t_1 = \frac 72$ and $t_2 = \frac 32$.
\end{example}

If each $e_k = 2$ in Theorem \ref{omnibus} and $m$ is as small as
possible, then we obtain an analogue of Sylvester's Theorem for
forms of even degree. 
\begin{corollary}\label{Sylv3can} 
\ 

(i) A general binary form of degree $d = 6s$ can be
written as
\begin{equation}\label{E:sylv6s}
\la x^{6s} + \sum_{j=1}^{2s}(\al_j x^2 + \be_j  x y + \ga_j y^2)^{3s}
\end{equation}
for some $\la \in \cc$.

(ii) A general binary form of degree $d = 6s+2$ can be
written as
\begin{equation}\label{E:sylv6s2}
\sum_{j=1}^{2s+1} (\al_j x^2 + \be_j  x y + \ga_j y^2)^{3s+1}.
\end{equation}

(iii) A general binary form of degree $d = 6s+4$ can be
written as
\begin{equation}\label{E:sylv6s4}
\la_1 x^{6s+4} + \la_2 y^{6s+4} +
\sum_{j=1}^{2s+1}(\al_j x^2 + \be_j  x y + \ga_j y^2)^{3s+2}
\end{equation}
for some $\la_i \in \cc$.
\end{corollary}

We have not been able to find an analogue to Sylvester's Algorithm
 for determining the representations \eqref{E:sylv6s},
 \eqref{E:sylv6s2}, \eqref{E:sylv6s4} in Corollary \ref{Sylv3can}. In
 the linear case,  $(\al x + \be y)^d$ is killed by $\be
 \frac{\partial}{\partial x} -  
\al \frac{\partial}{\partial y}$, and two operators of this shape
commute. Although  each  $(\al x^2 + 2\be xy + \ga y^2)^d$ is killed by
the non-constant-coefficient $(\be x + \ga y) \frac{\partial}{\partial x} - 
(\al x + \be y) \frac{\partial}{\partial y}$, two operators of this
kind do not usually commute. The smallest constant-coefficient
differential operator which kills $(\al x^2 + 2\be xy + \ga y^2)^d$
has degree $d+1$; the product of any two of these would kill every
form of degree $2d$ and so provide no information.

Let us say that \eqref{E:reveal} is a {\it neat} canonical form if $m=0$,
and of {\it Sylvester-type} if it is neat and if $e_k = e$ for $1 \le
k \le r$. Counting the numbers of neat and Sylvester-type canonical forms
leads to some  number theory. The first lemma is standard.

\begin{lemma}\label{egyptian}
Given  $0 < \frac pq \in \qq$ and $0 < n \in \nn$,
there exist only finitely many choices of $m_j \in \zz$, 
$0 < m_1 \le m_2 \cdots \le m_n$, such that $\frac pq = \sum_{j=1}^n
\frac 1{m_j}$.
\end{lemma}
\begin{proof}
If $n = 2$, then $\frac pq > \frac 1{m_1} \ge \frac p{2q}$ implies
that there are finitely many integral choices for $m_1$, each of which
determines $m_2 = (\frac pq - \frac 1{m_1})^{-1}$. 
Supposing the lemma valid for $n-1$, we have $\frac pq > \frac 1{m_1} \ge \frac 
p{nq}$, and each choice of $m_1$ implies the equation 
$\frac pq - \frac 1{m_1}= \sum_{j=2}^{n} \frac 1{m_j}$.
This has finitely many solutions by the induction hypothesis.
\end{proof}
\begin{theorem}\label{finitenice}
For fixed value of  $r$, there are only finitely many neat canonical
forms \eqref{E:reveal} with $r$ summands.
\end{theorem}

\begin{proof}
Suppose $m=0$ in Theorem \ref{omnibus}. Write $d = e_k m_k$, then by
\eqref{E:book}, 
\begin{equation}\label{E:neategyptian}
d+1 = \sum_{k=1}^r \left(\frac d{m_k} + 1\right) \implies
1 =  \sum_{k=1}^r \frac 1{m_{k}} +  \frac{r-1}d =
\sum_{k=1}^r \frac 1{m_{k}} + \sum_{\ell = 1}^{r-1} \frac 1d. 
\end{equation}
Now apply  Lemma \ref{egyptian} with $\frac pq = 1$ and $n=2r-1$:
there are only finitely many expressions of 1 as a sum of
$2r-1$ unit fractions, of which only a  subset satisfy the additional
restrictions of \eqref{E:neategyptian}.
\end{proof}

It is not hard to work out that for $r=2$, there are three neat canonical forms:
$(d,e_1,e_2) = (3,1,1)$, $(4,2,1)$ and $(6,3,2)$. The first is Theorem
\ref{Sylvcan}(i) with $d=3$, the second is Corollary \ref{6-22} with
$d=4$ (see Theorem \ref{6} below),
and the third is Theorem \ref{sextican}. 
When $r=3$, there are twenty-two neat canonical forms. 

Let $s(d)$ denote the number of neat Sylvester-type canonical forms of
degree $d$. Suppose $e_k = e$ for all $k$ in one of these. 
Then $e \ | \ d$ and, by \eqref{E:book}, $r(e+1) = d+1$, so
$(e+1) \ | \ (d+1)$.
Since $d \equiv 0 \pmod{e}$ and $d \equiv -1 \pmod{(e+1)}$, it follows
from the Chinese Remainder Theorem that $d \equiv e \pmod{e(e+1)}$;
that is,  $d=e + ue(e+1)$, $u \ge 1$, so that $e < \sqrt d$.
\begin{theorem}
 Let $S(N) := \sum_{d=1}^N s(d)$. Then $S(N)  = N + \mathcal O(N^{1/2})$
 and  $\sup_d s(d) = \infty$.
\end{theorem}
\begin{proof}
The generating function for the sequence $(s(d))$ is
\begin{equation}\label{E:414}
\sum_{n=1}^{\infty} s(d) x^d=
\sum_{e=1}^\infty\sum_{u=1}^\infty x^{e + ue(e+1)} = \sum_{e=1}^\infty
\frac{x^{e^2+2e}}{1-x^{e^2+e}} = \sum_{N=e}^\infty
\left\lfloor\frac{N-e}{e^2+e}\right\rfloor X^N.
\end{equation}
Let $T =\lfloor N^{1/2} \rfloor$.  It follows from \eqref{E:414} that
\begin{equation}\label{E:SNform}
S(N)= \sum_{n=1}^N s_n = 
\sum_{e=1}^\infty \left\lfloor \frac{N -e}{e^2+e} \right\rfloor =
\sum_{e=1}^T \left\lfloor \frac{N -e}{e^2+e} \right\rfloor. 
\end{equation}
Thus, using the telescoping sum for $\sum \frac 1{e(e+1)}$,
\eqref{E:SNform} implies that
\begin{equation}\label{E:416}
\begin{gathered}
S(N) \le \sum_{e=1}^T\frac{N - e}{e^2+e} = 
N\sum_{e=1}^T\frac 1{e^2+e} -
\sum_{e=1}^T\frac 1{e+1}
\\
\le N(1 - \tfrac 1{T+1}) - \log T + \mathcal O(1) = N - N^{1/2} + \mathcal
O(\log N).
\end{gathered}
\end{equation}
The lower bound is the same, minus $T$, so
\eqref{E:416} implies that $S(N) = N + \mathcal O(N^{1/2})$.

Now, $s(d)$ counts the number of $e<d$ so
that $e$ divides $d$ and $e+1$ divides $d+1$. If $d = 2^r-1$,
then $e+1 \ | \ 2^r$ implies that $e+1 = 2^t$ for some $t < r$. But
$2^t-1 \ |\ 2^r - 1$ if and only if $t \ | \ r$, hence $s(2^r-1) =
d(r) - 1$, where $d(n)$ denotes the divisor function. In particular,
$s(2^{2^t}-1) = t$, so the sequence $(s(d))$ is unbounded. More generally, if
$e \ | \ d$ and $e+1 \ | \ d+1$, then $e \ | \ d^2+2d$ and $e+1 \ | \
d^2+2d+1$, and since $e=d$ contributes to the count in $s(d^2+2d)$ but
not in $s(d)$, $s(d^2+2d) \ge s(d) + 1$.
\end{proof}
Half of the neat Sylvester forms come from 
Theorem \ref{Sylvcan}(i), another sixth come from Corollary
\ref{Sylv3can}(ii), etc. The smallest $d$ for which
$s(d)=2$ is $d=15$: $(e,r) = (1,8), (3,4)$, so a general binary form of
degree 15 is a sum of eight linear forms to the 15th power, or four cubics to
the 5th power. Mathematica computations show that the smallest $d$ for
which $s(d) = 3$   
 is $d=99$: $(e,r) = (1,50),(3,25),(9,10)$. For $d < 10^7$, the largest
 value of $s(d)$ is $s(7316000) = 12$. Note that $2^{2^{13}}- 1 =
 2^{4096}-1 \approx 1.04 \times 10^{1233}$, so the examples given in the proof
 are not likely to describe the fastest growth. We conjecture as well that
 $\{s(d)\}$ has an underlying distribution.

If the degree $d$ is prime, then Theorem \ref{e=1} accounts for all
canonical forms in Theorem \ref{omnibus}. 
The smallest $d$ which is not covered by Theorem \ref{e=1} is then $d=4$, and
there are two such cases, one of which is neat: $e_1=2, e_2=1, m=0$
and $e_1=2, m=2$. Both can be discussed constructively.

\begin{theorem}\label{6}
A general binary quartic  $p \in H_4(\cc^2)$ can be written as
\begin{equation}\label{E:2,1}
p(x,y) = (t_1 x^2 + t_2 x y + t_3 y^2)^2 + (t_4 x + t_5 y)^4
\end{equation}
in six different ways. Further, the set of possible values for
$\{\frac {t_5}{t_4}\}$ is the 
image of the set $\{0,\infty, 1,-1,i,-i\}$ under a M\"obius transformation.
\end{theorem}
\begin{proof}
By Theorem \ref{quarticgen}, if $p$ is a general binary
  quartic, then there exist $c_i, \la$ so that $p(c_1 x+ c_2 y,c_3 x + c_4
  y) =  p_{\la}(x,y) := x^4 + 6\la x^2y^2 + y^4$.
If \eqref{E:2,1} holds for $p_{\la}$, then
\begin{equation}
\begin{gathered}\label{E:swamp}
1 = t_1^2 + t_4^4,\quad  0 = 2t_1t_2 + 4t_4^3t_5,\quad  6\la = 2t_1t_3 +t_2^2+ 6
t_4^2t_5^2,\\ 0 = 2t_2t_3 + 4t_4t_5^3,\quad  1 = t_3^2+ t_5^4.
\end{gathered}
\end{equation}
First suppose that $t_4=0$. Then \eqref{E:swamp} implies that $1 =
t_1^2$ and  $0 = 2t_1t_2$, so $t_1 = 1$ (without loss of generality)
and $t_2 = 0$. The remaining equations
imply that $t_3 = 3\la$ and $t_5^4 = 1 - 9\la^2$. A similar
argument works if $t_5 = 0$, giving two representations:
\begin{equation}\label{E:easy}
p_{\la}(x,y) = (x^2 + 3\la y^2)^2 + (1-9\la^2)y^4 = (3\la x^2 + y^2)^2
+ (1-9\la^2)x^4. 
\end{equation}

Now suppose $t_4t_5 \neq 0$, so $t_1t_2t_3 \neq 0$ and
so 
\begin{equation*}
\frac{t_3}{t_1} = \frac{-2t_2t_3}{-2t_1t_2}
=\frac{4t_4t_5^3}{4t_4^3t_5} = \frac{t_5^2}{t_4^2} \implies 
\frac{1-t_3^2}{1-t_1^2}=\frac{t_5^4}{t_4^4}=
\frac{t_3^2}{t_1^2}\implies t_1^2 = t_3^2.
\end{equation*}
It follows that  $t_5 = i^k t_4$ and $t_3 = (-1)^kt_1$, and
\eqref{E:swamp} can be completely solved: 
\begin{equation*}
\begin{gathered}
t_4^4 = 1-t_1^2,\quad t_2 = 2i^k(t_1-t_1^{-1}),\quad 2 + 6(-1)^k\la = 4t_1^{-2}.
\end{gathered}
\end{equation*}
After some massaging of the algebra, this gives four representations:
\begin{equation}\label{E:base}
\begin{gathered}
p_{\la}(x,y) = \left(\frac {(-1)^k 2}{3\la+(-1)^k}\right)\left(x^2 -
  i^{3k}(3\la - (-1)^k) x y + 
 (-1)^k y^2\right)^2 \\ + 
\left(\frac{3\la - (-1)^k}{3\la + (-1)^k}\right) \left( x + i^ky
\right)^4, \quad k = 0,1,2,3.  
\end{gathered}
\end{equation}
In order to find the six representations of $p$ as \eqref{E:2,1}, we
start with the six representations of $p_{\la}$ given in
\eqref{E:easy} and \eqref{E:base}, in which $t_4x 
+ t_5y$ is a multiple of one of the six linear forms
$x,y,x+i^ky$. Apply the the inverse of the map
$(x,y) \mapsto (c_1 x+ c_2 y,c_3 x + c_4 y)$, which takes $t_4 x + t_5
y$ to a multiple of
$t_4(c_4 x - c_2 y) + t_5(-c_3 x + c_ 1 y)$: $\frac {t_5}{t_4}
\mapsto G(\frac {t_5}{t_4})$, where $G(z) = \frac{c_1z - c_2}{c_4 -
  c_3 z}$.
\end{proof}

\begin{theorem}\label{4;2}
Given two fixed non-proportional binary linear forms $\ell_1, \ell_2$,
a general binary 
quartic in $H_4(\cc^2)$ has two representations as 
\begin{equation}\label{E:4,1,1}
p(x,y) = (t_1 x^2 + t_2 x y + t_3 y^2)^2 + t_4 \ell_1(x,y)^4 + t_5 \ell_2(x,y)^4.
\end{equation}
\end{theorem}
\begin{proof}
Given $p, \ell_1,\ell_2$, make an invertible linear change of variable
taking $(\ell_1,\ell_2) \mapsto (x,y)$, and suppose $p(x,y) \mapsto
q(x,y)= \sum_ia_ix^{4-i}y^i$. Then $q$ has the shape 
\eqref{E:4,1,1} if and only if the coefficients of $x^3y, x^2y^2,x
y^3$ in $(t_1 x^2 + t_2 x y + t_3 y^2)^2$ and $q$ agree. Thus, we seek
to solve the system 
\begin{equation}\label{E:grub}
a_1 = 2t_1t_2,\quad a_2 = 2t_1t_3+t_2^2, \quad a_3 = 2t_2t_3.
\end{equation}
But \eqref{E:grub} implies $a_1t_2^2  - 2a_2 t_1t_2 +  2a_3
t_1^2 = 0$, hence in general, there are exactly two values of $\be$ so
that $t_2 = \be t_1$; in each case, $t_1^2 = \frac{a_1}{2\be}$. The two choices
of sign for $t_1$ lead to the same square, and $t_3 = \frac{a_1}{a_3} t_1$, so
\eqref{E:grub} has these two solutions.
\end{proof}

\begin{table}

\begin{center}

\begin{tabular}{c c c c c}

$d$ & $e_1,\dots,e_r$ & m & $F(d;e)$ & Source \\

\hline

\hline

$d$ & $1^{\lfloor \frac{d+1}2 \rfloor}$ & 0 or 1 &1 & Theorem \ref{Sylvcan} \\
$d$ & $1^r$ & $d+1-2r$ & 1 & Theorem \ref{e=1}\\

4 & 2,1 & 0 & 6 & Theorem \ref{6} \\

4 & 2 & 2  & 2 & Theorem \ref{4;2} \\

6 & 3,2 & 0 & 40 & \cite{V1,V2} \\

6 & 2,$1^2$ & 0 & 22 & Experiment \\

6 & 3,1 & 1 & 14 & Experiment \\

6 & $2^2$ & 1 & 9 & Experiment \\

6 & 2,1 & 2 & 12 & Experiment \\

6 & 3 & 3 & 5 & Experiment \\

6 & 2 & 4 & 5 & Experiment \\

8 & 2,$1^3$ & 0 & 62 & Experiment \\

10 & 2,$1^4$ & 0 &  147 & Experiment \\

12 & 2,$1^5$ & 0 &  308 & Experiment \\

$2s$ & 2,$1^{s-1}$ &0 &  $2\binom{s+3}5 + \binom{s+2}3$ & Conjecture

\\

\\

\end{tabular}
\caption{Values of $F(d;e)$}

\end{center}

\end{table}

In the case of Theorem \ref{omnibus} let $F(d;e_1,\dots,e_r)$ denote
the number of different representations that a general $p \in
H_d(\cc^2)$ has, by our convention. We present in Table 1 a complete list of 
proved or conjectural values when $d \le
6$, reflecting numerical experiments on Mathematica.
 (Recall that if $d$ is prime, then Theorem \ref{e=1} presents all
possible canonical forms of this type.) The conjectural value of
$F(2s;2,1^{s-1})$ is suggested by the given data for $2 
\le s \le 6$ and OEIS\cite[A081282]{O}.

V\'arilly-Alvarado, in \cite{V1,V2}, constructs explicitly all 240
representations of $x^6 + y^6$ as $f^2 + g^3$; he considers forms
multiplied by roots of unity as different, which explains the
appearance of $\frac {240}{2\cdot 3}$ in the table above. This
is also proved to be the number of representations for a general sextic.

To describe the experiments for  $F(2s;2,1^{s-1})$ more precisely, we
generate a form  
\begin{equation*}
p(x,y)= \sum_{k=0}^{2s}\binom {2s}k a_k x^{2s-k}y^k,
\end{equation*} 
where $a_k = t + iu$ for random integers $t,u$ in $[-100,100]$.
In case $m \le 2$, we assume a change of variables so that the fixed
linear forms are $x^d$ or $y^d$; for $m>2$ we choose
additional linear forms with random coefficients. Let $h(x,y) = U x^2
+ V x y + W y^2$ for variables $(U,V,W)$ and 
let  $q(x,y) = p(x,y) - h^s(x,y)$,  and apply Sylvester's Algorithm to
$q$. That is, we construct the 
$(s+2) \times s$ matrix $A_{s-1}(q)$, with polynomial entries in 
 $(U,V,W)$ of degree $s$ and require that it have
rank $< s$.  This is done by counting the number of
$(U,V,W)$ which are common zeros of all $s \times s$ minors. This
number is divided by $s$ to account for $h^s = (\zeta_s^k h)^s$.
As a back of the envelope calculation, one might take the first $s-1$
rows of $A_{s-1}$ and use the cofactors to compute a non-trivial
null-vector. Ignoring possible cancellation, the components would be
polynomials of degree $s(s-1)$ in 
$(U,V,W)$. Taking the dot product with the last three rows of
$A_{s-1}$ gives three polynomials of degree $s^2$. Ignoring
cancellations and multiplicity, there should be 
$(s^2)^3$ common zeros, and dividing by $s$ gives an upper
bound for $F(2s;2,1^{s-1})$ of $s^5$. The conjectural value
is asymptotically $\frac 1{60}s^5$, which shows the same order
of growth.

\section{Quadratic forms and the proof of Theorem 1.2}
We begin this section with a constructive proof of  Theorem
\ref{uppertri} which will serve as a template for constructive proofs
involving cubic 
forms. 
\begin{proof}[Constructive Proof of Theorem \ref{uppertri}]
Suppose $p \in H_2(\cc^n)$, and specifically,
\begin{equation*}
p(x_1,\dots,x_n) =\sum_{i=1}^n a_{ii}x_i^2 + 2\sum_{1 \le i<j \le n} a_{ij}x_ix_j.
\end{equation*}
Then $\frac{\partial p}{\partial x_1} = 2\sum_{j=1}^n a_{1j}x_j$,
Since $a_{11} \neq 0$  in general, we can define
$q(x_1,\dots,x_n) = p(x_1,\dots,x_n) - \frac 1{a_{11}}(
  \sum_{j=1}^n a_{1j}x_j )^2$.
Observe that $\frac{\partial q}{\partial x_1} = 0$,
so  $q = q(x_2,\dots,x_n)$.
Iterating this argument gives the construction. There is only one
linear form $\pm \ell$ so that $\frac{\partial p}{\partial x_1}
=2\ell\frac{\partial \ell}{\partial x_1}$, so
the representation is unique. 
\end{proof} 

Constant-counting for sums of squares is complicated by the action of
the orthogonal group on a sum of $t$ squares. If $M \in Mat_t(\cc)$
and $MM^t = I$, then 
\begin{equation*}
\sum_{i=1}^t f_i^2 = \sum_{i=1}^t\left(\sum_{j=1}^t m_{ij}f_j\right)^2.
\end{equation*} 
When $t=2$, choose $\theta \in \cc$ and let
 $e^{i\theta} = \cos\theta + i \sin \theta := (u,v)$, so that
 \begin{equation}\label{E:rotate}
f^2 + g^2 = (uf - vg)^2 + (vf +ug)^2.
\end{equation}
This means that we may safely remove one monomial from one of the summands.
\begin{theorem}\label{so2s}
A general binary form $p \in H_{2s}(\cc^2)$ can be written as
\begin{equation}\label{E:so2s}
 \left( \sum_{k=0}^s t_kx^{s-k}y^k \right)^2 +  \left(
  \sum_{k=1}^s t_{s+k}x^{s-k}y^k \right)^2.  
\end{equation}
in $\binom{2s-1}s$ different ways.
\end{theorem}
\begin{proof}
The non-constructive proof is a simple
application of Corollary \ref{jake}. Writing
\eqref{E:so2s} as $f^2 + g^2$ gives the partials with respect to the $t_j$'s as 
\begin{equation*}
\left \{ 2x^{s-k}y^k f,\  0 \le k \le s \right \} \cup 
\left \{ 2x^{s-k}y^k g,\  1 \le k \le s \right \}; 
\end{equation*}
specializing to $f = x^s$ and $g = y^s$ above gives all monomials in
$H_{2s}(\cc^2)$.  

The more obvious expression
\begin{equation}\label{E:bigso2s}
p(x,y) = f^2(x,y) + g^2(x,y),\qquad g, h \in H_s(\cc^2)
\end{equation}
is not a canonical form, because $2(s+1) > 2s+1$. 
However, every  sum of two squares can be formally factored, and these behave
nicely with respect to \eqref{E:rotate}.
\begin{equation*}
\begin{gathered}
f^2+g^2 = (f+ig)(f-ig) \iff  \\ (uf + vg)^2 + (vf -ug)^2 = \left(
  e^{i\ta}(f + ig)\right) \left(  e^{-i\ta}(f - ig)\right). 
\end{gathered}
\end{equation*}
Suppose $p(1,0) = a_0 \neq 0$ (true for general $p$) and
\eqref{E:bigso2s} holds, where $f(1,0) = \rho$ and  $g(1,0) =
\tau$. Then $\rho^2 + \tau^2 = a_0$, so that $\frac{\tau}{\rho} \neq \pm
i$ and the coefficient of $x^s$ in
$vf+ug$ will be $v\rho + u\tau = \sin\ta \rho + \cos \ta \tau$, which
is zero exactly when $\tan \ta =  -\frac{\tau}{\rho}$. Thus for
precisely one value of $\tan \ta$, the right-hand side of
\eqref{E:rotate} will be in the form  \eqref{E:so2s}. This
determines a pair $(\pm u,\pm v)$; however, the squares in \eqref{E:so2s}
will be the same.

In other words, each distinct factorization of $p$ (up to multiple) as a product
of two $s$-ic forms, when combined with the orthogonal action of
\eqref{E:rotate}, 
yields exactly one representation as \eqref{E:so2s}. A general $p \in
H_{2s}(\cc^2)$ is a product of $2s$ distinct linear factors; these can
be organized  into an unordered pair of products of $s$ distinct
linear factors in $\frac 12\binom {2s}s = \binom{2s-1}s$ ways.  
\end{proof}
The ``lost" degree of freedom in a sum of squares
never arises in Theorem \ref{omnibus} because $2(\frac d2 + 1) >
d+1$. The missing monomial $x^s$ in the second summand of
\eqref{E:so2s} may be replaced by any specified monomial $x^{s-k_0}y^{k_0}$ by
a similar argument.

Another classical result is that a general ternary quartic is a sum of three
squares of quadratic forms, generally in 63 different ways up to the
action of the orthogonal group (see \cite{PRSS}.)
Hilbert proved that {\it every} positive
semidefinite  $p \in H_4(\rr^3)$ is a sum of three squares from
$H_2(\rr^3)$ \cite{H}. He then showed that there exist psd forms in
$H_6(\rr^3)$ and $H_4(\rr^4)$ which are not sums of squares in
$H_3(\rr^3)$ and $H_2(\rr^4)$, respectively, which ultimately led to
his 17th problem. (See \cite{Re1} for much more on this subject.)

 A constructive discussion of Hilbert's theorem on  $p \in H_4(\rr^3)$
 has recently been given in papers  
by Powers and the author \cite{PR00}, Powers, Scheiderer, Sottile and the author
\cite{PRSS}, Pfister and Scheiderer \cite{PS} and Plaumann, Sturmfels
and Vinzant \cite{PSV}. A non-constructive proof (without the count)
can easily be given. 
\begin{theorem}\label{so3s}
A general ternary quartic $p \in H_4(\cc^3)$ can be written as $p =
q_1^2 + q_2^2 + q_3^2$, where $q_j \in H_2(\cc^3)$. 
\end{theorem}
\begin{proof}
We take $q_i$'s so that the monomial $x^2$ only appears in $q_1$ and
the monomial $y^2$ only appears in $q_1$ and $q_2$, and so the number of
coefficients in the $q_j$'s is 
$6+5+4=15$. Taking the partials where $(q_1,q_2,q_3) = (x^2,y^2,z^2)$
shows that $J$ contains $2x^2\{x^2,y^2,z^2,xy,xz,yz\}$,
$2y^2\{y^2,z^2,xy,xz,yz\}$ and $2z^2\{z^2,xy,xz,yz\}$, and so is
equal to $H_4(\cc^3)$. 
\end{proof}
Since $3\binom{m+2}2-3 < \binom{2m+1}2$ for $m \ge 3$,  this
result does not generalize to ternary forms of higher even degree. 

The situation is somewhat simpler over $\rr$. 
A real version of Theorem \ref{so2s} appears in \cite{R1}. If $p$ is
real and positive definite and $p = f^2 + g^2$, where $f$ and $g$ are
also real, then the factors of $p$ consist of $s$ conjugate pairs. In the
factorization $p = (f+ig)(f-ig)$, the pairs must be split between
the conjugate factors, and if $p$ has distinct factors, this can be
done in $2^{s-1}$ different ways.  
A real generalization of Theorem \ref{so3s} appears in \cite[Corollary
2.12]{CLR}.
Suppose a real psd form $p \in
H_{2s}(\rr^n)$ is a sum of $t$ squares and $x^{\be_i} \in H_s(\rr^n),
1 \le i \le t$, is given. Then there is a representation $p =
\sum_{j=1}^t g_j^2$, in which $x^{\be_i}$ does not occur in $g_j$ for $j
> i$. This argument can also be applied to a {\it general} sum of $t$
squares over $\cc$, but it no longer applies to all forms. For
example, if $xy = (ax + by)^2 + (cx + dy)^2$, then $abcd \neq 0$.

\section{Cubic forms and proofs of Theorems 1.3 and 1.4}
In this section, we present three representations for forms in
$H_3(\cc^n)$ as a sum of cubes of linear forms. The first two are
canonical; the third isn't, but it represents {\it all} cubics, not
just general cubics.

We begin with Theorem \ref{reichcan}, which first  appeared \cite{Re1}
in a 1987 paper of Boris Reichstein. At the time of
this writing, \cite{Re1}  has had no citations in MathSciNet. (It was
discussed in \cite{R5} and, from there, in \cite{CM}. The former was
never submitted for publication and the latter appeared in an unindexed
journal.) The original presentation and proof in \cite{Re1} were given
for trilinear forms (see $\S2$); the theorem is applied to cubic forms
there mainly in the examples. 

By iterating \eqref{E:reich}, we obtain a canonical form 
for $p \in H_3(\cc^n)$, see \cite[p.98]{Re1}.
\begin{corollary}\label{Reichstein2}
A general $n$-ary cubic $p \in H_3(\cc^n)$ can be written uniquely as
\begin{equation}\label{E:fullreich}
p(x_1,\dots,x_n) = \sum_{m=0}^{\lfloor (n-1)/2\rfloor}
\sum_{k=1}^{n-2m}(t^{\{k\}}_{m,1+2m}x_{1+2m} + \cdots + t^{\{k\}}_{m,n} x_n)^3 
\end{equation}
for some $t^{\{k\}}_{m,j} \in \cc$.
\end{corollary}
This gives $p$ as a sum of $n + (n-2) + \cdots  =
\bigl\lfloor\frac{(n+1)^2}4\bigr\rfloor$
cubes. 
Recall that by Alexander-Hirschowitz, for $n \neq 5$, a general cubic
form in $n$ variables can be written as a sum of $ \bigl \lceil
  \frac{(n+1)(n+2)}6 \bigr\rceil$ cubes. Thus  \eqref{E:fullreich}
is a canonical form which  represents a general 
cubic as a sum of about 50\% more cubes than the true minimum; this is
due to the large number of linear forms with restricted sets of
variables. 

Reichstein's proof of Theorem \ref{reichcan} requires the well-known
``generalized eigenvalue problem'' for pairs of symmetric matrices, as
interpreted for quadratic forms: 
if a general pair of quadratic forms $f,g \in H_2(\mathbb C^n)$ is
given, then there exist  $n$ linearly independent forms
$L_i(x) = \sum_{j=1}^n \al_{ij}x_j$ and $c_i \in \mathbb C$ so that
\begin{equation}\label{E:simult}
f = \sum_{i=1}^n L_i^2, \qquad g = \sum_{i=1}^n c_i L_i^2.
\end{equation}
If $M_f, M_g$ are the matrices associated to $f,g$, then the $c_i$'s
are the $n$ roots of the determinantal equation $\det (M_g - \la M_f)
= 0$, which are generally distinct, so the $L_i$'s are uniquely
determined up to multiple. We may also assume that the coefficients
$\al_{ij}$ of the linear forms are generally non-zero; cf. Corollary
\ref{reichugly}.

\begin{proof}[Proof of Theorem \ref{reichcan}]
For general $p \in H_3(\mathbb C^n)$, we
 simultaneously diagonalize
 $f = \frac{\partial p}{\partial x_1}$ and  $g=\frac{\partial p}{\partial
   x_2}$ as in \eqref{E:simult}. Since mixed partials are equal,
\begin{equation}\label{E:64}
\frac{\partial f}{\partial x_2} = 
\frac{\partial g}{\partial x_1} = 
\sum_{i=1}^n 2\al_{i2}L_i = \sum_{i=1}^n 2c_i\al_{i1}L_i,
\end{equation}
and since the $L_i$'s are linearly independent, \eqref{E:64} implies
that $\al_{i2} = c_i\al_{i1}$. 

It is generally true that $\al_{i1} \neq 0$. Let
\begin{equation*}
q(x_1,\dots,x_n) = p(x_1,\dots,x_n) - \sum_{i=1}^n \tfrac
1{3\al_{i1}} L_i^3.  
\end{equation*}
It follows that
\begin{equation*}
\begin{gathered}
\frac{\partial q}{\partial x_{1}} = \frac{\partial p}{\partial x_{1}}
-  \sum_{i=1}^n \frac{3\al_{i1}}{3\al_{i1}}L_i^2 = 
 \frac{\partial p}{\partial x_{1}} -  \sum_{i=1}^n L_i^2 = 0 , \\ 
\frac{\partial q}{\partial x_{2}} = \frac{\partial p}{\partial x_{2}}
- \sum_{i=1}^n \frac{3\al_{i2}}{3\al_{i1}} L_i^2= \frac{\partial p}{\partial x_{2}}
- \sum_{i=1}^n c_i L_i^2 = 0.
\end{gathered}
\end{equation*}
Since $\frac{\partial q}{\partial x_{1}} =\frac{\partial q}{\partial
  x_{2}} = 0$, we have $q = q(x_3,\dots,x_n)$. 

For uniqueness, suppose \eqref{E:reich} holds and
$\ell_k(x_1,\dots,x_n) = \sum_j\be_{kj}x_j$. Then 
\begin{equation*}
 \frac{\partial p}{\partial x_{1}} = \sum_{k=1}^n 3\be_{k1}
 \ell_k^2;\qquad  
 \frac{\partial p}{\partial x_{2}} = \sum_{k=1}^n 3\be_{k2} \ell_k^2. 
\end{equation*} 
Thus, after a scaling, $\frac{\partial p}{\partial x_{1}}$ and
$\frac{\partial p}{\partial x_{2}}$ have already been simultaneously
diagonalized (as in \eqref{E:simult}), and the $\ell_k$'s are, up to
multiples, a rearrangement of the $L_k$'s. 
\end{proof}
We now give a constructive proof of Theorem \ref{slinkycan}, which
gives a different canonical form for $H_3(\cc^n)$ requiring even more
cubes.
\begin{proof}[Proof of Theorem \ref{slinkycan}]
The constant-counting makes this a potential
canonical form: the variables are $t_{\{i,j\},k}$ with $1 \le i \le j
\le k \le n$, and there are $\binom{n+2}3 = N(n,3)$ such triples
$(i,j,k)$. Given $p \in H_3(\cc^n)$, $\frac{\partial p}{\partial x_n}$ is a
quadratic form, so we can 
generally complete the square by Theorem \ref{uppertri}: 
\begin{equation*}
\frac{\partial p}{\partial x_n} = \sum_{j=1}^n (t_{jj}x_j + \cdots +
t_{jn}x_n)^2. 
\end{equation*}
Then $t_{jn} \neq 0$ for general $p$ and if we let 
\begin{equation*}
q(x_1,\dots,x_n) = p(x_1,\dots,x_n) - \sum_{j=1}^n  \tfrac
1{3t_{jn}} (t_{jj}x_j  + \cdots + t_{jn}x_n)^3,
\end{equation*}
then $\frac{\partial q}{\partial x_n} = 0$, so  $q = q(x_1,\dots,x_{n-1})$.
Iterate this construction to get \eqref{E:slinky}.

Uniqueness follows by working backwards. If \eqref{E:slinky} holds for
a cubic $p$, then it gives $\frac{\partial p}{\partial x_n}$ in its
(unique) upper-triangular diagonalization. This can be integrated with
respect to $x_n$ and subtracted from $p$, giving a cubic
$q(x_1,\dots,x_{n-1})$. Again, iterate.
\end{proof}
It is not hard to give nonconstructive proofs of Theorems
\ref{reichcan} and \ref{slinkycan} using Corollary
\ref{jake}. These are left for the reader.

We first presented this next construction in \cite{R5}; an outline of
the proof can be found in \cite{CM}. This is not a canonical form, but
is included here because it gives an absolute upper bound for the
length of cubic forms.

\begin{theorem}\label{slowpoke}
If $p \in H_3(\cc^n)$, then there exists an invertible linear change of
variables $y_j = \sum \la_{jk} x_k$ and $n$ linear forms $\ell_j$ so that
for some $q \in H_3(\cc^{n-1})$,
\begin{equation}\label{E:slowpoke}
p(x_1,\dots,x_n) = \sum_{j=1}^n \ell_j^3(x_1,\dots,x_n) + q(y_2,\dots,y_n).
\end{equation}
Thus every cubic in $n$ variables is a sum of at most $\binom {n+1}2$
cubes of linear forms.
\end{theorem}
\begin{proof}
Define linear forms $\ell_{j,m}(y)$ for $1 \le j \le m+1$ by 
\begin{equation}\label{E:predrab}
\begin{gathered}
\ell_{j,m}(y_1,\dots,y_n)= y_j + \al \sum_{j=1}^m y_j, \qquad 1 \le j \le m, \\
\ell_{m+1,m}(y_1,\dots,y_n)=-(1+m\al) \sum_{j=1}^m y_j, \qquad \al =
\frac{-(m+1)+\sqrt{m+1}}{m(m+1)}. 
\end{gathered}
\end{equation}
Then it can be easily checked that
\begin{equation}\label{E:drab}
\sum_{j=1}^{m+1} \ell_{j,m}(y) = 0 \quad \text{and}  \quad
\sum_{j=1}^{m+1} \ell_{j,m}^2(y) = 
\sum_{k=1}^m y_k^2.
\end{equation}

Suppose $0 \neq p \in H_3(\mathbb C^n)$. Use Biermann's Theorem
to find a point $u$ where $p(u) \neq 0$, and after an
invertible linear change of variables, taking $\{x_j\} \mapsto
\{u_j\}$, we may assume that $p(1,0,\dots,0) = 1$ and so 
\begin{equation}\label{E:tedious}
p = u_1^3 + 3h_1(u_2,\dots,u_n)u_1^2 + 3h_2(u_2,\dots,u_n)u_1 +
h_3(u_2,\dots,u_n),
\end{equation}
where $deg\ h_j = j$.
Now let $u_1 = y_1 - h_1(u_2,\dots, u_n)$ to clear the quadratic
term, so
\begin{equation}
p = y_1^3 + 3y_1\tilde h_2(u_2,\dots,u_n) +
\tilde h_3(u_2,\dots,u_n), 
\end{equation}
where again $deg\ \tilde h_j = j$.
Diagonalize  $\tilde h_2(u_2,\dots,u_n)$ as a quadratic form into $y_2^2 +
\dots + y_r^2$, where $r \le n$, and make the accompanying change of
variables. We now have 
\begin{equation}
p = y_1^3 + 3y_1(y_2^2 + \dots + y_r^2) + k_3(y_2,\dots, y_n);\quad r \le n,
\end{equation}
where $deg\ k_3 = 3$.
Finally, using \eqref{E:predrab} and \eqref{E:drab}, we construct $g$,
a sum of $r \le n$ cubes: 
\begin{equation}
\begin{gathered}
g(y_1,\dots,y_n) := \frac 1{r} \sum_{j=1}^r \left(y_1 + \sqrt{r}\cdot
  \ell_{j,r-1}(y_2,\dots,y_r)\right)^3 \\
= \frac 1r \sum_{j=1}^r y_1^3 + \frac 3{\sqrt r}  \sum_{j=1}^r y_1^2\ell_{j,r-1} +
3\sum_{j=1}^r y_1\ell^2_{j,r-1} + \sqrt{r}  \sum_{j=1}^r\ell^3_{j,r-1}
\\
=y_1^3 + 3y_1(y_2^2 + \dots + y_r^2) + \sqrt{r}
\sum_{j=1}^r\ell^3_{j,r-1}(y_2,\dots,y_r). 
\end{gathered}
\end{equation}
Then $q:=p - g$ is a cubic form in $(y_2,\dots,y_n)$ as in
\eqref{E:slowpoke}.  Iteration of this argument
shows that any cubic $p \in H_3(\cc^n)$ is a sum of at most
$\frac{n(n+1)}2$ cubes. 
\end{proof}

Theorem \ref{E:reich} can be extended to a 
canonical form for quartics as a sum of fourth powers of
linear forms. Note that $x_n$ appears in each summand of \eqref{E:fullreich}, 
with, generally, a non-zero coefficient. 

\begin{corollary}\label{reichugly}
For general $p \in H_4(\cc^n)$, there exist $\ell_k 
 \in H_1(\cc^n)$ and $q \in  H_4(\cc^{n-1})$ so that, with $a(n) =
 \lfloor \frac{(n+1)^2}4 \rfloor$, 
\begin{equation*}
p(x_1,\dots,x_n) = \sum_{k=1}^{a(n)} \ell_k(x_1,\dots,x_n)^4 +
q(x_1,\dots,x_{n-1}). 
\end{equation*}
As a consequence, a general $p \in H_4(\cc^n)$ can be written as
\begin{equation*}\label{E:fullreichugly}
p(x_1,\dots,x_n) = \sum_{m=0}^{\lfloor (n-1)/2\rfloor}\sum_{r=1+2m}^n
\sum_{k=1}^{r-2m}(t^{\{k\}}_{m,r,1+2m}x_{1+2m} + \cdots + t^{\{k\}}_{m,r,r} x_r)^4. 
\end{equation*}
\end{corollary}
\begin{proof}
By Corollary \ref{reichcan} and \eqref{E:fullreich}, for general $p
\in H_4(\cc^n)$, we can write 
\begin{equation}\label{E:reichbonus}
\begin{gathered}
\frac{\partial p} {\partial x_n} =  \sum_{m=0}^{\lfloor (n-1)/2\rfloor}
\sum_{k=1}^{n-2m}(t^{\{k\}}_{m,1+2m}x_{1+2m} + \cdots +
t^{\{k\}}_{m,n} x_n)^3 \\ =: \sum_{m=0}^{\lfloor (n-1)/2\rfloor}
\sum_{k=1}^{n-2m} (\ell_m^{\{k\}}(x))^2.
\end{gathered}
\end{equation}
As before, if $q = p - \sum_{k,m} \frac 1{4t^{\{k\}}_{m,n}}
\ell_{k,m}^4$, then
$\frac{\partial q} {\partial x_n} = 0$, so  
$q =  q(x_1,\dots,x_{n-1})$. Repeat as before.
There are $N(n,3)$ coefficients in \eqref{E:reichbonus}, and since
$N(n,3)+N(n-1,4)=N(n,4)$, the count is correct for a canonical form.  
\end{proof}
Note that there is no variable which appears in each linear form
in \eqref{E:reichbonus}, so the argument can't be extended to
quintics. For the same reason, Theorem \ref{slinkycan} does not
extend to quartics. By combining Theorems \ref{reichcan} and
\ref{reichugly}, we obtain canonical
forms as a sum of powers of linear forms in the four exceptional
cases of Alexander-Hirschowitz, of course at the expense of the number
of summands. With regards to ternary quartics and Theorem
\ref{notclebsch}, Corollary \ref{reichugly} becomes  the
following canonical form for $H_4(\cc^3)$ as a sum of seven fourth powers. 
\begin{equation*}
\begin{gathered}
\sum_{k=1}^3 (t_{k1}x_1 + t_{k2}x_2+t_{k3}x_3)^4 + t_{10} x_3^4 +
\sum_{\ell = 1}^2 (u_{\ell 1}x_1 +  
u_{\ell 2}x_2)^4 + u_5 x_1^4.
\end{gathered}
\end{equation*}

There is an arithmetic obstruction to a ``Reichstein-type'' canonical form
for quartics; that is, one in which each linear form is allowed to
involve each variable.  If 
\begin{equation}\label{E:reichgen}
p(x_1,\dots,x_n) = 
\sum_{k=1}^r(\al_{k1}x_1 + \cdots + \al_{kn} x_n)^4 + q(x_1,\dots,x_m).
\end{equation}
were a canonical form for some $n$, then we would have  $N(n,4) = rn +
N(m,4)$.  However, for $n = 12$, there does not exist $m < 12$ so that
$12 \ | \ \binom{15}4 - \binom{m+3}4$, so no such canonical form can
exist. More generally, let 
\begin{equation}
A_d = \left\{n :  0 \le m < n \implies n \not |\  \tbinom{n+d-1}d -
  \tbinom{m+d-1}d\right\}
\end{equation}
denote the set of $n$ for which this argument rules out
Reichstein-type canonical forms.
We present without proof a number of results about $A_d$. Note that
there is no obstacle for \eqref{E:reichgen} in prime degree, such as $d=2,3$.
\begin{proposition}\label{noreich}
\

(i) If $3 \not | \ k$, then $n = 2^{2k}\cdot 3 \in A_4$.

(ii) If $p \equiv 1 \pmod{144}$ is prime, then $12p \in A_4$.

(iii) If $p$ is prime, then $p\ |\ \tbinom{n+p-1}p -  \tbinom{n}p$, hence
  $A_p = \emptyset$ for prime $p$ .
  
(iv)  The smallest elements of $A_6, A_8, A_{10}, A_{12}, A_{14}$ and
  $A_{15}$ are 10, 1792, 6, 242, 338 and 273 respectively. If $A_9$ or
  $A_{16}$ are non-empty, then their smallest elements are at least
  $10^5$. 
\end{proposition}

\section{Subspace canonical forms and the Proof of Theorem 1.11}

One natural generalization of the definition of canonical forms is
to consider maps $F: X \mapsto H_d(\cc^n)$ where $X \subset \cc^M$ is an
$N(n,d)$-dimensional subspace of $\cc^M$. (Similar ideas can be found
in Wakeford \cite{W}, though his approach is different from ours.) 
These can be analyzed
 in the simplest non-trivial case: $M=4, N(2,2)=3$.

\begin{proof}[Proof of Theorem \ref{hyperplane}]
Assume that some $c_j \neq 0$. Without loss of generality, we may
assume that $c_4 
\neq 0$ and divide through by $c_4$ so that the equation is $t_4 =
a_1t_1 + a_2 t_2 + a_3 t_3$, where $a_i = -c_i/c_4$ for $i = 1,2,3$.
Then \eqref{E:4} can be  parameterized as a map from $\cc^3\mapsto
H_2(\cc^2)$ as: 
\begin{equation}\label{E:cand}
F(t;x) = (t_1 x + t_2 y)^2 + (t_3 x +(a_1 t_1 + a_2 t_2 + a_3 t_3) y)^2.
\end{equation}
The partials with respect to the $t_j$'s are:
\begin{equation}\label{E:mire}
\begin{gathered}
2x(t_1 x + t_2 y) + 2a_1y (t_3 x +(a_1 t_1 + a_2 t_2 + a_3 t_3) y), \\
2y(t_1x + t_2 y) + 2a_2y (t_3 x +(a_1 t_1 + a_2 t_2 + a_3 t_3) y), \\
2(x + a_3 y)(t_3 x +(a_1 t_1 + a_2 t_2 + a_3 t_3) y).
\end{gathered}
\end{equation}
Now, \eqref{E:cand} is a canonical form if and only if there exists 
a choice of $t_i$ so that the three quadratics in \eqref{E:mire}
span $H_2(\cc^2)$. A computation shows that the determinant of the
forms in \eqref{E:mire} 
 with respect to the basis $\{x^2,xy,y^2\}$ is the cubic
\begin{equation}\label{E:muck}
-8 ((a_1 a_2 - a_3) t_1 + (1 + a_2^2)t_2  + (a_2 a_3+a_1) t_3) (a_1 t_1^2 + 
   a_2 t_1 t_2 + a_3 t_1 t_3 - t_2 t_3).
\end{equation} 
The second factor in \eqref{E:muck} always has the term $-t_2t_3$ and
so never vanishes, 
hence this determinant is not  identically zero (and \eqref{E:cand} is
a canonical form), unless 
\begin{equation}\label{E:badplane}
a_1a_2 - a_3 = 1 + a_2^2 = a_2a_3 + a_1 = 0.
\end{equation}
In the exceptional case where \eqref{E:badplane} holds, then $a_2 = \ep$, where $\ep = \pm i$,
and $a_3 = \ep a_1$. Evaluating \eqref{E:cand} at $(x,y) = (a_1,\ep)$
yields
\begin{equation*}
\begin{gathered}
(a_1t_1 + \ep t_2)^2 + (a_1t_3 + \ep a_1 t_1 + \ep^2t_2 + \ep^2a_1t_3)^2
 \\ = (a_1t_1 + \ep t_2)^2 + ((1+\ep^2)a_1t_3 + \ep a_1 t_1 + \ep^2t_2
 )^2 = (a_1t_1 + \ep t_2)^2 + \ep^2  ( a_1 t_1 + \ep t_2)^2 = 0,
\end{gathered}
\end{equation*}
as claimed. 
\end{proof}
It would be interesting to know how Theorem \ref{hyperplane}
generalizes to higher degrees. 

Conjecture \ref{zerosum} is true for degree 2 by Theorem
\ref{hyperplane}. We have verified it for even degrees up to eight 
by Corollary \ref{jake} applied to random choices for $\al_j, \be_j$
in \eqref{E:sylvalt2}. We  hold some hope that generalizations such as
Conjecture \ref{zerosum} will have applications in more than two
variables as well.


 \end{document}